\documentclass[11pt]{article}

\usepackage{amssymb,amsmath,latexsym}
\usepackage[dvips]{graphics}
\usepackage{mathrsfs}
\usepackage{float}

\newtheorem{thm}{Theorem}[section]
\newtheorem{prop}{Proposition}[section]

\newtheorem{lemma}{Lemma}[section]
\numberwithin{equation}{section}
\renewcommand{\baselinestretch}{1.4}
\def\ZZ{\mathbb{Z}}
\def\ii{\mathsf{i}}
\def\HH{\mathbb{H}}

\def\NN{\mathbb{N}}
\def\EE{\mathbb{E}}
\def\LLL{\mathscr{L}}

\def\RR{{\bf R}}

\def\FF{{\mathscr{F}}}
\def\GG{\mathscr{G}}
\def\AA{\mathscr{A}}
\def\RR{{\bf R}}

\def\MM{{\textrm{M}}}

\def\ang #1{{\langle #1\rangle}}
\def\<{\langle}
\def\>{\rangle}
\def\pf{\noindent{\bf Proof.} }

\def\E{{\bf E}}

\def\P{{\bf P}}

\def\qed{{\hfill $\Box$\medskip}}

\allowdisplaybreaks[4]

\input{epsf.sty}
\usepackage{epsfig}

\begin{document}
\title{\bf The Complete Convergence Theorem Holds for Contact Processes in a Random Environment on $\ZZ^d\times
\ZZ^+$}

\author{Qiang Yao\footnote{School of Finance and Statistics, East China Normal University, Shanghai 200241, China;
                      E-mail: qyao@sfs.ecnu.edu.cn.}
       ~and Xinxing Chen\footnote{Department of Mathematics, Shanghai Jiaotong University, Shanghai 200240, China.}
        }

\maketitle{}

\begin{abstract}
In this article, we consider the basic contact process in a static random environment on the half space $\ZZ^d\times
\ZZ^+$ where the recovery rates are constants and the infection rates are independent and identically distributed random variables. We show that, for almost every environment, the complete convergence
theorem holds. This is a generalization of
the known result for the classical contact process in the half space case.

 \end{abstract} \noindent{\bf 2000 MR subject
classification:} 60K35


\section{Introduction}

The aim of this paper is to obtain the complete convergence theorem
for the contact process in a random environment on the half space $(\HH,\EE)$. The vertex set is
$\HH=\ZZ^d\times\ZZ^+~(d\geq1)$, where $\ZZ=\{0,\pm 1,\pm
2,\cdots\}$ denotes the set of integers and $\ZZ^+=\{0,1,2,\cdots\}$
denotes the set of nonnegative integers. And the edge set is $\EE=\{(x,y):~x,y\in\HH,~\|x-y\|=1\}$, where $\|\cdot\|$ denotes the
Euclidean norm. Here, we treat the graph as unoriented; that is, $(x,y)$ and $(y,x)$ denote the same edge for all $x,y\in\HH$ satisfying $\|x-y\|=1$. The environment is given by $\lambda=(\lambda_e)_{e\in\EE}$, a collection of nonnegative random variables
which are indexed by the edges in $\EE$. The random
variable $\lambda_e$ gives the infection rate on edge $e$. We let the law of
$(\lambda_e)_{e\in\EE}$ be independent and identically distributed with law $\mu$, which puts mass $1$ on $[0,+\infty)$. To
describe the environment more formally, we consider the following
probability space. We take
$\Omega_1=[0,+\infty)^{\EE}$ as the sample space, whose elements are represented by
$\omega=(\omega(e):~e\in\EE)$. The value $\omega(e)$ corresponds to
the infection rate on edge $e$; that is,
$\lambda_e(\omega)=\omega(e)$ for every $e\in\EE$. We take $\FF_1$
to be the $\sigma$-field of subsets of $\Omega_1$ generated by the
finite-dimensional cylinders. Finally, we take product measure on
$(\Omega_1,\FF_1)$; this is the measure
$\P^{\mu}=\prod\limits_{e\in\EE}\mu_e$, where $\mu_e$ is a measure
on $[0,+\infty)$ satisfying $\mu_e(\omega(e)\in\cdot)=\mu(\cdot)$
for every $e\in\EE$. The probability space
$(\Omega_1,\FF_1,\P^{\mu})$ describes the environment.



Next, we fix the environment $\lambda=(\lambda_e)_{e\in\EE}$ and consider the basic contact process under this environment. The state space of the contact process
$\xi=\xi(\lambda)$ is $\{A:~A\subseteq\HH\}$, and the transition rates are as
follows:
$$\left\{\begin{array}{ll} \xi_{t}\rightarrow\xi_{t}\setminus\{x\} \text{ for } x\in\xi_{t} \text{ at rate } 1,\\
\xi_{t}\rightarrow\xi_{t}\cup\{x\} \text{ for } x\notin\xi_{t}
\text{ at rate } \sum_{y: \|y-x\|=1} \lambda_{(y,x)}\textbf{1}_{\{y\in
\xi_t\}}.\end{array}\right.$$  Readers can refer to the standard references Liggett \cite{Liggett1985} and Durrett \cite{Durrett1988} for how these rates rigorously determine a Markov process $\xi(\lambda)$ on $(\Omega_2,\FF_2,\P_{\lambda})$ and for much on the contact process as well as other interacting particle systems. Denote by $\xi^A(\lambda)$ the process with initial state $A$.
If $\lambda$ is random, then the transition rates
are random variables and therefore $\P_{\lambda}$ becomes a random measure. We say
that $\xi^A$ survives if $\xi^A_t\neq \emptyset$ for all $t\ge 0$,
while $\xi^A$ dies out if there exists $t>0$ such that
$\xi^A_t=\emptyset$.

The model in several special environments have been studied before. For example, Bezuidenhout and Grimmett \cite{Bezuidenhout-Grimmett1990} studied the case when $\mu(\{c\})=1$ for some $c>0$. (In fact, this is an almost nonrandom environment.) Bramson et al. \cite{Bramson-Durrett-Schonmann1991} studied the case when $\mu(\{a,b\})=1$ for some $0<a<b$. Chen and Yao \cite{Chen-Yao2009} studied the case when $\mu(\{0,c\})=1$ for some $c>0$. All the above models belong to static environments; that is, the environment does not change as time goes. There are some models concerning contact processes in dynamic environments; see, for example, Broman \cite{Broman2007}, Remenik \cite{Remenik2008}, and Steif and
Warfheimer \cite{Steif-Warfheimer2008}.\\

Regarding complete convergence, Bezuidenhout and Grimmett \cite{Bezuidenhout-Grimmett1990} showed that the complete convergence theorem holds for the
basic contact process on $\ZZ^d$. Chen and Yao \cite{Chen-Yao2009} showed that the complete convergence theorem holds for the contact process on
open clusters of half space $\ZZ^d\times\ZZ^+$. In this paper, we will show that, for the general model described above, the complete convergence theorem still holds for almost every environment. It generalizes the results of Bezuidenhout and Grimmett \cite{Bezuidenhout-Grimmett1990} and  Chen and Yao \cite{Chen-Yao2009} in the half space case. Denote by $\nu_\lambda$ the upper invariant measure,
that is, the weak limit of the distribution of
$\xi_t^{\HH}(\lambda)$ as $t\rightarrow\infty$, and denote by
$\delta_{\emptyset}$ the probability measure which puts mass one on
the empty set. Note that, since $\lambda$ is random, $\nu_{\lambda}$ is a random measure. We then have the following complete convergence theorem, which is the main result of this paper.
\begin{thm}\label{t:1.1}  Suppose $\mu$ puts mass $1$ on
$[0,\infty)$. Then there exists $\Omega_0\subseteq\Omega_1$ with
$\P^{\mu}(\Omega_0)=1$, such that for all $\omega\in\Omega_0$ and
$A\subseteq \HH$,
$$
\xi_t^A(\lambda)\Rightarrow
\nu_{\lambda}\cdot\P_{\lambda}(\xi^A(\lambda)
{\rm~~survives~})+\delta_{\emptyset}\cdot\P_\lambda(\xi^A(\lambda)~~{\rm
dies~~ out~})
$$
as $t$ tends to infinity, where \lq~$\Rightarrow$\rq ~stands for
$\P_{\lambda}$-weak convergence.\\
\end{thm}

The main purpose of this paper is to prove Theorem \ref{t:1.1}, which will be specified in the following sections. The rest of this paper is organized as follows. In Section 2, we give some preliminaries including some basic notation, together with an introduction to the important \lq graphical representation\rq. In Section 3, we prove the \lq block conditions\rq~ which are essential to the proof of Theorem \ref{t:1.1}. We prove it under three different cases. In Section 4, we use these blocks to construct the route and use the renormalization method to make further preparations. Finally, in Section 5, we prove Theorem \ref{t:1.1} by checking the two equivalent conditions in Theorem 1.12 of \cite{Liggett1999}.

The main idea of the whole procedure is enlightened by Bezuidenhout and Grimmett \cite{Bezuidenhout-Grimmett1990}. But there are some big differences. In order to make good use of some symmetric properties, we need to consider the annealed law first (Sections 3 and 4), then go back to the quenched law to get the desired result (Section 5). The fact is, under the annealed law, the process is not Markovian, but events depending on disjoint subgraphs are relatively independent. In consequence, we can only get \lq space blocks\rq~ rather than \lq space-time blocks\rq~ as in Bezuidenhout and Grimmett \cite{Bezuidenhout-Grimmett1990}. Furthermore, we can only use these \lq space blocks\rq~ to obtain the result in the half space case. We believe that the result will hold for the whole space case, but we cannot construct the independent \lq restart process\rq~ as in Bezuidenhout and Grimmett \cite{Bezuidenhout-Grimmett1990} by adopting the method of this paper.

\section{Preliminaries}
We only prove the case $d=1$; that is, $\HH=\ZZ\times\ZZ^+$. Our technique still
works for the case $d\ge 2$ after trivial modifications. In this section, we introduce some basic notation for the following analysis.

When $d=1$, for simplicity we use a complex number $a+b\ii$ to denote the vertex $(a,b)\in\HH=\ZZ\times\ZZ^+$, where $a\in\ZZ$ and $b\in\ZZ^+$. Furthermore, we use the notation $\lceil a+b\ii,c+d\ii\rfloor$ to denote the rectangle $$[\min\{a,c\},\max\{a,c\}]\times[\min\{b,d\},\max\{b,d\}],$$ that is, $a+b\ii$ and $c+d\ii$ are diagonal sites of this rectangle. The notation $\lceil\cdot\rfloor$ can be used in a more flexible way. If $a=c$~(respectively, $b=d$) then $\lceil a+b\ii,c+d\ii\rfloor$ denotes a vertical~(respectively, horizontal) line. We can also let $a$, $b$, $c$, or $d$ be infinity. For example, $\lceil-3,3+\infty\ii \rfloor$ denotes the infinite \lq rectangle\rq~$[-3,3]\times[0,+\infty)$.

Now, we introduce a special
notation  $\ang{\cdot,\cdot}$. For $a,c\in\ZZ$ and $b,d\in\ZZ^+$, define
$$\ang{a+b\ii,c+d\ii}:=\left\{\begin{array}{ll} \left\{(u,v)\in \EE: u,v\in \lceil a+b\ii,c+d\ii\rfloor, \{\Re(u),\Re(v)\}\not \subseteq\{a,c\}\right\},~~~~\; \text{ if } |a-c|\ge 2|b-d|,\\
\left\{(u,v)\in \EE: u,v\in \lceil
a+b\ii,c+d\ii\rfloor,\{\Im(u),\Im(v)\}\not \subseteq\{b,d\}\right\},
~~~~\; \text{ if } 2|a-c|\le |b-d|.\end{array}\right.$$ Then
$\ang{a+b\ii,c+d\ii}$ is an edge set. See Figure 1.
\begin{figure}[H]\label{picc:1}
\center
\includegraphics[width=10.0true cm, height=2.9true cm]{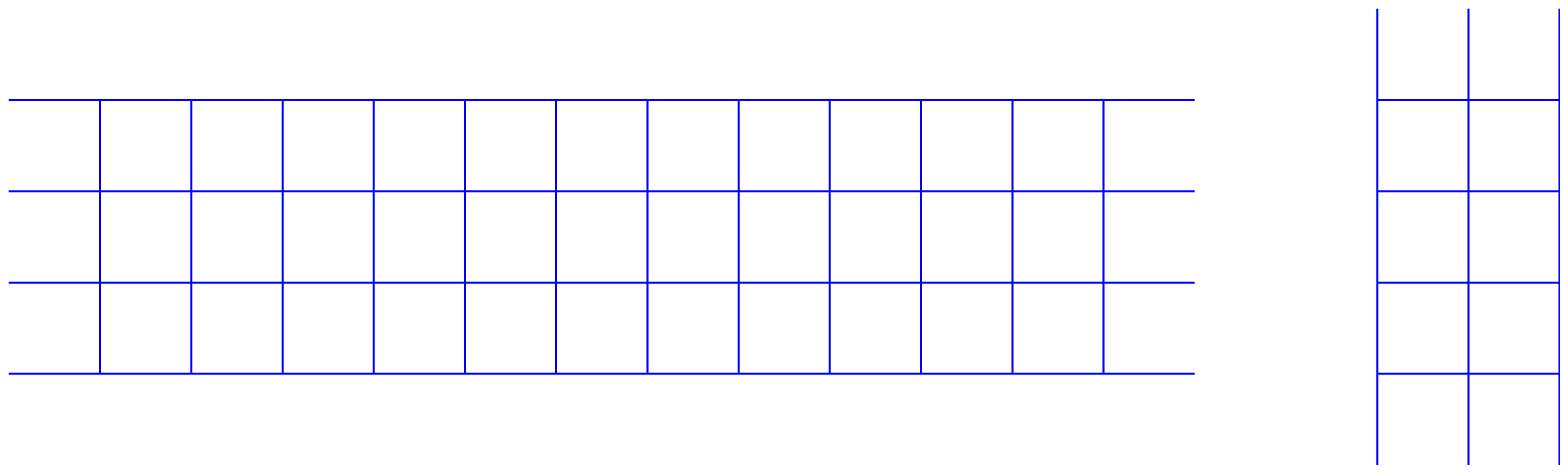}
\caption{$\ang{a+b\ii,c+d\ii}$ }
\end{figure}

For a real number $a$, let $[a]$ be the largest integer which is no larger than $a$. Then for $x\in \HH$ and $M\in \ZZ^+$, set $$B_x(M):=\lceil x-M-M\ii, x+M+M\ii\rfloor\cap\HH$$ to be the \lq ball\rq~centered at $x$ and with radius $M$ (but restricted on $\HH$).

Denote by $\P$ a probability measure which satisfies
$$
\P\left(\xi^A\in
\cdot\right)=\int\P_{\lambda}\left(\xi^A(\lambda)\in
\cdot\right)\P^\mu(d\omega).
$$
We call $\P$ the annealed (average) law and $\P_{\lambda}$ the quenched law. Note that the contact process is Markovian under the quenched law, while it is not Markovian under the annealed law.\\

We shall make abundant use of the graphical representation of the
contact process which was first proposed in Harris \cite{Harris1978}. We follow the
notation of Bezuidenhout and
Grimmett \cite{Bezuidenhout-Grimmett1990}. Fix $\lambda$, and think of
the process as being imbedded in space-time. Along each \lq
time-line\rq~ $x\times [0,\infty)$ are positioned \lq deaths\rq~at
the points of a Poisson process with intensity 1, and between each
ordered pair $x_1\times[0,\infty)$, $x_2\times[0,\infty)$ of
adjacent time-lines are positioned edges directed from the first to
the second having centers forming a Poisson processes of intensity
$\lambda_{(x_1,x_2)}$ on the set
$\frac{1}{2}(x_1+x_2)\times[0,\infty)$.  These Poisson processes are
taken to be independent of each other. The random graph obtained
from $\HH\times[0,\infty)$ by deleting all points at which a death
occurs and adding in all directed edges can be used as a percolation
superstructure on which a realization of the contact process is
built. We shall make free use of the language of percolation. For
example, for $A,B\subseteq \HH\times [0,\infty)$, we say that $A$ is
joined to $B$ if there exists $a\in A$ and $b\in B$ such that there
exists a path from $a$ to $b$ traversing time-lines in the direction
of increasing time (but crossing no death) and directed edges
between such lines; for $C\subseteq\HH\times[0,\infty)$, we say that
$A$ is joined to $B$ within $C$ if such a path exists using segments
of time-lines lying entirely in $C$.
We next extend the notion \lq within\rq~ in this paper. For
 $A,B\subseteq\HH\times[0,\infty)$ and $C\subseteq\HH$, we say that $A$ is joined to $B$
within $C$ if such a path exists using segments of time-lines lying
entirely in $C\times[0,\infty)$; for $D\subseteq \EE$, we say that $A$ is
joined to $B$ within $D$ if such a path exists using directed edges
having centers lying entirely in $D'\times[0,\infty)$, where
$D'=\{\frac{x_1+x_2}{2}: (x_1,x_2)\in D\}$.

 For $x\in\HH, r\in \ZZ^+$ and $t\in
[0,\infty)$,
  we call $(x\times t)_r$ a horizontal (respectively, vertical) seed with $2r+1$ sites if all sites
in $\lceil x-r,x+r\rfloor$ (respectively, $\lceil x-r\ii, x+r\ii\rfloor$)
are infected at time $t$. We say that a horizontal seed $(x\times s)_r$
is joined to a vertical seed $(y\times t)_r$ if $\lceil
x-r,x+r\rfloor\times s$ is joined to
 $z\times t$ for all $z\in \lceil y-r\ii,y+r\ii\rfloor$.  The word
 \lq seed\rq~ comes from Grimmett  \cite{Grimmett1999}.

\section{Block conditions}
To prove Theorem \ref{t:1.1}, we need to get the \lq block conditions\rq~
for the survival of the process. The construction is enlightened by
Bezuidenhout and Grimmett \cite{Bezuidenhout-Grimmett1990}, and was
used successfully in the proof of the complete convergence theorem
for contact processes on open clusters of $\ZZ^d\times\ZZ^+$; see
Chen and Yao \cite{Chen-Yao2009}. We first introduce some notation we will need.

For
$h,w\in\NN$, define the random set
$$ \Phi^R(h,w):=\{x\in \lceil w, w+h\ii\rfloor:
\lceil-r,r\rfloor\times 0\text{~ is~ joined~ to~} x\times [0,\infty)
\text {~within~}\lceil-w,w+h\ii\rfloor \},
$$
 Hence,
$\Phi^R(h,w)$ is a subset of the right side of the box
$\lceil-w,w+h\ii\rfloor$. Similarly, define $ \Phi^L(h,w)$ as a subset
of the left side. Define the random set $\Phi^{UR}(h,w)$, which is a subset of
the right part of the up side, as follows:
$$ \Phi^{UR}(h,w):=\{x\in \lceil h\ii,w+h\ii\rfloor: \lceil-r,r\rfloor\times 0\text{~ is~ joined~ to~} x\times
[0,\infty) \text {~within~}\lceil-w,w+h\ii\rfloor\times[0,\infty)
\}.
$$
 Similarly, define
$\Phi^{UL}(h,w)$ as a subset of the left part. Furthermore, denote
\begin{equation}\label{e:3.11}\Phi(h,w):=\Phi^L(h,w)\cup\Phi^R(h,w)\cup\Phi^{UL}(h,w)\cup\Phi^{UR}(h,w).\end{equation}
Then we have
\begin{equation}\label{e:4.5}
|\Phi(h,w)|\le|\Phi^{L}(h,w)|+|\Phi^{UL}(h,w)|+
|\Phi^{UR}(h,w)|+|\Phi^{R}(h,w)|\le |\Phi(h,w)|+3.\\
\end{equation}

Next, we present the \lq block conditions\rq~ in the following proposition.
\begin{prop}\label{p:3.2}
Suppose that $\P(\xi^0~\text{survives})>0$. Then, for any $N\in\NN$ and $\varepsilon>0$ sufficiently small, one of the
following two
assertions must be true.\\
$(\verb"1")$ There exist constants $h,w$ with $w=4h$, such that
\begin{equation}\label{e:3.2}
\P(|\Phi^R(h,w)|>N)>1-\varepsilon,~~~\P(|\Phi^R(h,
2w)|>N)>1-\varepsilon.
\end{equation}
$(\verb"2")$ There exist constants $h , w $ with $~ 8h\geq w$, such
that
\begin{equation}\label{e:3.3}
\P(|\Phi^{UR}(h,w)|>N)>1-\varepsilon,~~~ \P(|\Phi^R(2h,
w)|>N)>1-\varepsilon.
\end{equation}
Here, $|\cdot|$ denotes the cardinality of a set.\\
\end{prop}

The content of Proposition \ref{p:3.2} is quite similar to Lemma 3.2 in Chen and Yao \cite{Chen-Yao2009}, but things are much more difficult here. In the Bernoulli bond percolation model, it is easy to get the property that the existence of crossing from bottom to top of a box is small if the ratio of the height to the width of the box is large enough. However, in the model presented in this
paper, this property is not obvious. So we need to develop some new
ideas to make the construction. In detail, we consider the following three
cases, which will be proved in Sections 3.1--3.3, respectively. Here and henceforth, for any $A,B\subseteq\HH$ we say that $\xi^A$ survives within $B$ if, for any $t>0$, there exists $x\in B$ such that $A\times0$ is joined to $x\times t$ within $B$, while we say that $\xi^A$ dies out within $B$ otherwise.

\textbf{Case 1.} $\mu(\{0\})>0$.

\textbf{Case 2.} $\mu(\{0\})=0$ and $\xi^0$ cannot survive within any
\lq slab\rq~ $\lceil-k, k+\infty\ii\rfloor$ with positive
probability.

\textbf{Case 3.} $\xi^0$ survives within some \lq slab\rq~ with positive probability.\\

The following lemma is important to the analysis throughout this paper. The idea of its proof comes from the Remark on page 347 of  \cite{Steif-Warfheimer2008}.

\begin{lemma}\label{l:3.1} If $\P(\xi^0~\text{survives})>0$, then
\begin{equation}\label{e:3.1} \lim\limits_{r\rightarrow\infty}\P(\xi^{\lceil-r,r\rfloor}~\text{survives})=1.
\end{equation}
\end{lemma}
\pf Let $Y_x:=\textbf{1}_{\{\xi^x~\text{survives}\}}$ for any $x\in(-\infty,+\infty)$. Then, by our assumption, we have $$\P(Y_x=1)=\P(\xi^0~\text{survives})$$ for any $x\in(-\infty,+\infty)$. Furthermore, it follows from the graphical representation that $\{Y_x\}_{x\in(-\infty,+\infty)}$ is ergodic. So
$$\P(\xi^{\lceil-r,r\rfloor}~\text{survives})=\P(\exists x\in\lceil-r,r\rfloor~\text{s.t.}~\xi^x~\text{survives})\rightarrow\P(\exists x\in(-\infty,+\infty)~\text{s.t.}~Y_x=1)=1$$ as $r$ tends to infinity, as desired.\qed

\subsection{Proof of Case 1}
In this subsection, we shall prove that the block conditions hold if $\mu(\{0\})>0$. By Lemma \ref{l:3.1}, for any $\varepsilon>0$ sufficiently small we can take some
$r\in \NN$  such that
\begin{equation}\label{e:4.3}
\P(\xi^{\lceil-r,
r\rfloor}\text{~survives})>1-\frac{\varepsilon^6}{4}.
\end{equation}
Set $w_n=2^n$ and $h_n= 2^{w_n^2}$ for each $n>100r$. 
Since $\mu(\{0\})>0$, we have that, for
sufficiently large $n$, with large probability there exists $1<h<h_n-1$
such that $\lambda_{(x,x+\ii)}=0$ for all $x\in \lceil
-w_n+h\ii,w_n+h\ii\rfloor$. Obviously, if $\lambda_{(x,x+\ii)}=0$
for all $x\in \lceil -w_n+h\ii,w_n+h\ii\rfloor$, then
$$\Phi^{UL}(h_n,w_n)=\Phi^{UR}(h_n,w_n)=\emptyset.$$ So we can conclude that there exists $n_0$ such that, for $n>n_0$,
\begin{equation}\label{e:4.4}\P(|\Phi^{UL}(h_n,w_n)|+|\Phi^{UR}(h_n,w_n)|=0)>1-\frac{\varepsilon^6}{2}.
\end{equation}
Let $\FF_n$ denote the $\sigma$-field generated by the graphical representation within $\lceil-w_n,w_n+h_n\ii\rfloor~(n=1,2,\cdots)$. Note that, for any $n\in\NN$, if $\lambda_e=0$ for all $ e\in\{(x,y):~x\in \Phi(h_n,w_n),
y\not\in\lceil-w_n,w_n+h_n\ii\rfloor\} $, then
$\xi^{\lceil-r,r\rfloor}$  must die out, since
 no  sites outside  $\lceil-w_n, w_n+h_n\ii\rfloor$ can be infected.
This implies that
$$
\P\left(\left.\xi^{\lceil-r,r\rfloor} \text{ dies
out}~\right|~\FF_n\right)\geq [\mu(\{0\})]^{|\Phi(h_n,w_n)|+2}
$$
for any $n\in\NN$. By the martingale convergence theorem,
$$\P\left(\left.\xi^{\lceil-r,r\rfloor} \text{ dies
out}~\right|~\FF_n\right)\rightarrow\textbf{1}_{\{\xi^{\lceil-r,r\rfloor} \text{ dies
out}\}}~~\text{a.s.}$$ as $n$ tends to infinity. Since $0<\mu(\{0\})<1$, it follows that
$$\lim\limits_{n\rightarrow\infty}|\Phi(h_n,w_n)|=\infty~\text{almost surely on}~\{\xi^{\lceil-r,r\rfloor}\text{ survives}\}.$$
Therefore,
$$
\P(\exists m, \forall n>m,|\Phi(h_n,w_n)|>
2N~|~\xi^{\lceil-r,r\rfloor}\text{ survives} )=1.
$$
Hence there exists $n_1>n_0$ such that, for $n>n_1$,
\begin{equation}\label{e:4.6}
\P(|\Phi(h_n,w_n)|> 2N~|~\xi^{\lceil-r,r\rfloor}\text{
survives})>1-\frac{\varepsilon^6}{4}.
\end{equation}
By (\ref{e:4.3}) and (\ref{e:4.6}), if $n>n_1$, then
\begin{equation}\label{e:4.61}
\P(|\Phi(h_n,w_n)|>2N)>1-\frac{\varepsilon^6}{2}.
\end{equation}
Furthermore, from (\ref{e:4.5}) we can see that $|\Phi(h_n,w_n)|>2N$ and $|\Phi^{UL}(h_n,w_n)|+|\Phi^{UR}(h_n,w_n)|=0$ together imply that $|\Phi^L(h_n,w_n)|+|\Phi^R(h_n,w_n)|>2N$. Therefore, by (\ref{e:4.4}) and (\ref{e:4.61}), we get that, if $n>n_1$, then
\begin{align*}
&\P(|\Phi^L(h_n,w_n)|+|\Phi^R(h_n,w_n)|>2N)\\
\geq&\P(|\Phi(h_n,w_n)|>2N,~|\Phi^{UL}(h_n,w_n)|+|\Phi^{UR}(h_n,w_n)|=0)\\
\geq&\P(|\Phi(h_n,w_n)|>2N)+\P(|\Phi^{UL}(h_n,w_n)|+|\Phi^{UR}(h_n,w_n)|=0)-1\\
>&1-\varepsilon^6.
\end{align*}
Using the Fortuin--Kasteleyn--Ginibre (FKG) inequality~(see Theorem 2.4 of
Grimmett \cite{Grimmett1999}) and the symmetric property, we can get
\begin{align*}
\varepsilon^6&>\P(|\Phi^L(h_n,w_n)|+|\Phi^R(h_n,w_n)|\leq 2N)\\
&\geq\P(|\Phi^L(h_n,w_n)|\leq N , |\Phi^R(h_n,w_n)|\leq N)\\
&\geq[\P(|\Phi^R(h_n,w_n)|\leq N)]^2
\end{align*}
Consequently, when $n$ is large,
\begin{equation}\label{e:4.7}
\P(|\Phi^R(h_n,w_n)|> N)>1-\varepsilon^3~(>1-\varepsilon^2>1-\varepsilon).
\end{equation}
Similarly, we have
\begin{equation}\label{e:4.71}
\P(|\Phi^R(h_n,2w_n)|> N)>1-\varepsilon^3~(>1-\varepsilon^2>1-\varepsilon)
\end{equation}
when $n$ is large.

Comparing (\ref{e:4.7}) and (\ref{e:4.71}) with (\ref{e:3.2}), we see that the ratio of $h_n$ to $w_n$ is
much larger than we want. Hence we need to reduce the height. Let $k'_n=w_n^2-n+2$ and $h'_n=h_{n}/2^{k'_{n}}$ for
$n=1,2,\cdots$. Then $4 h_n'=w_n$. If
$$\P(|\Phi^R(h'_n,w_n)|>N)>1-\varepsilon^2,~~~\P(|\Phi^R(h'_n,
2w_n)|>N)>1-\varepsilon^2$$
 for some
$n$, then $(\verb"1")$ is true. Otherwise, at least one of the two following
statements
must be true.

$(\verb"3")$ There exists a subsequence $(n_i)$ such that $
\P(|\Phi^R(h_{n_i}',w_{n_i})|>N)\leq 1-\varepsilon^2$.

$(\verb"4")$ There exists a subsequence $(n_i)$ such that $
\P(|\Phi^R(h_{n_i}',2w_{n_i})|>N)\leq 1-\varepsilon^2$.

For $i=1,2,\cdots$, take $w_{n_i}'=w_{n_i}$ if $(\verb"3")$ is true, and take $w_{n_i}'=2w_{n_i}$ if $(\verb"4")$ is true. Then, for any $i$, we have
$$w_{n_i}'\leq 8h_{n_i}'~~\text{and}~~\P(|\Phi^R(h_{n_i}',w_{n_i}')|>N)\leq 1-\varepsilon^2.$$
Meanwhile, from (\ref{e:4.7}) and (\ref{e:4.71}), we get
$$\P(|\Phi^R(h_{n_i},w_{n_i}')|> N)>1-\varepsilon^2$$ for any $i$. So, for any $i$, there exists $0\leq k\leq k'_{n_i}$ such that
$$\P\left(\left|\Phi^R\left(\frac{h_{n_i}}{2^{k+1}},w_{n_i}'\right)\right|> N\right)\leq1-\varepsilon^2,~~\P\left(\left|\Phi^R\left(\frac{h_{n_i}}{2^k},w_{n_i}'\right)\right|> N\right)>1-\varepsilon^2.$$
Set $h_i^*=h_{n_i}/2^{k+1}$ and $w_i^*=w_{n_i}'$. It follows that
\begin{equation}\label{e:4.8}
w_i^*\leq8h_i^*,~~\P(|\Phi^R(2h_i^*,w_i^*)|> N)>1-\varepsilon^2~\text{and}~
\P(|\Phi^R(h_i^*,w_i^*)|> N)\leq1-\varepsilon^2
\end{equation}
for any $i$.

We next show that there exists $i_0$ such that
\begin{equation}\label{e:4.9}\P(|\Phi^{UL}(h_{i_0}^*,w_{i_0}^*)|+|\Phi^{UR}(h_{i_0}^*,w_{i_0}^*)|>
2N)>1-\varepsilon^2.
\end{equation} In fact, if no such $i_0$ exists, then
 $\P(|\Phi^{UL}(h_i^*,w_i^*)|+|\Phi^{UR}(h_i^*,w_i^*)|> 2N)\leq 1-\varepsilon^2$ for all $i$. Using (\ref{e:4.5}), (\ref{e:4.8}), and the FKG
 inequality, we can get that, for any $i$,
\begin{align*}
&\P(|\Phi(h_i^*,w_i^*)|\leq
4 N-3) \\
\geq &\P(|\Phi^L(h_i^*,w_i^*)|\leq N,~|\Phi^R(h_i^*,w_i^*)|\leq N,~|\Phi^{UL}(h_i^*,w_i^*)|+|\Phi^{UR}(h_i^*,w_i^*)|\leq2N )\\
\geq&\P(|\Phi^L(h_i^*,w_i^*)|\leq N)\cdot\P(|\Phi^R(h_i^*,w_i^*)|\leq
N)\cdot\P(|\Phi^{UL}(h_i^*,w_i^*)|+|\Phi^{UR}(h_i^*,w_i^*)|\leq2 N)\\
\geq& \varepsilon^6.
\end{align*}
However, $h_i^*$ tends to infinity as $i\rightarrow\infty$. This
implies that there exists a strictly increasing subsequence
($h_{i_j}^*$) such that
\begin{equation}\label{e:4.81}
\P(|\Phi(h_{i_j}^*,w_{i_j}^*)|\leq 4N-3)\geq \varepsilon^6.
\end{equation}
On the other hand, by an argument similar to that of (\ref{e:4.6}), we have that, when $j$ is sufficiently large,
\begin{equation}\label{e:4.82}
\P(|\Phi(h_{i_j}^*,w_{i_j}^*)|>4N-3~|~\xi^{\lceil-r,r\rfloor}\text{ survives})>1-\frac{3\varepsilon^6}{4}.
\end{equation}
(\ref{e:4.3}) and (\ref{e:4.82}) together imply that, when $j$ is sufficiently large,
\begin{equation}\label{e:4.83}
\P(|\Phi(h_{i_j}^*,w_{i_j}^*)|>4N-3)>1-\varepsilon^6.
\end{equation}
(\ref{e:4.83}) contradicts (\ref{e:4.81}). As a result, (\ref{e:4.9}) is true for some $i_0$.

Let $h^*=h_{i_0}^*, w^*=w_{i_0}^*$. Then (\ref{e:4.9}) together with the FKG inequality and the symmetric property lead to
$$
\P(|\Phi^{UL}(h^*,w^*)|> N)=\P(|\Phi^{UR}(h^*,w^*)|>
N)>1-\varepsilon.
$$
So $(\verb"2")$ is true, and the proof of Case 1 is completed.\qed

\subsection{Proof of Case 2}

%

In this subsection we shall prove that the block conditions hold if $\mu(\{0\})=0$ and if $\xi^0$ cannot survive within any
\lq slab\rq~ $\lceil-k, k+\infty\ii\rfloor$ with positive
probability. Fix $N\in\NN$ and $\varepsilon>0$ sufficiently small. By Lemma \ref{l:3.1}, we can take some
$r\in \NN$  such that
\begin{equation}\label{e:5.3}
\P(\xi^{\lceil-r,
r\rfloor}\text{~survives})>1-\frac{\varepsilon^2}{16}.
\end{equation}
Set
\begin{equation}\label{e:5.31}
E:=\P( 0\times 0\text{~~is~joined~to~}
z\times1\text{~~within~}\{0\}\cup\lceil 1,4N+N\ii\rfloor\text{~~for~
all~} z\in \lceil 4N,4N+N\ii\rfloor)
\end{equation}
and $$\alpha:=\P(E);$$ then $\alpha>0$. Let $U$ be large enough to
ensure that, in $[U/20N]$ or more independent trials of an experiment with
success probability $\alpha $, the probability of obtaining at least
one success exceeds $1-\frac{\varepsilon}{4}$.
 Let $a$ be the minimal value which satisfies
$\mu([a,\infty))>1-\frac{\varepsilon}{200N^2}$. Then, for any set
$A\subset \EE$ with $\#A\le 20 N^2$,
\begin{equation}\label{e:5.32}
\P^\mu(\lambda_e\ge a, e\in
A)=(\mu([a,\infty)))^{\#A}>\left(1-\frac{\varepsilon}{200N^2}\right)^{20N^2}>1-\frac{\varepsilon}{8}.
\end{equation}
The value of $a$ is strictly larger than
$0$, since $\mu((0,\infty))=1$. Set
$$
\beta:=\P(0\times 0 \text{~~is ~joined~
to~} z\times 1 \text{~~within~} \{0\}\cup\lceil1, 4N+N\ii\rfloor
\text{~~for~ all~}z\in \lceil 4N,4N+N\ii\rfloor~|~\lambda_e=a\text{ for all }e\in\EE).
$$
Then $\beta>0$, since $a>0$. Let $V$ be large enough to ensure that, in $[V/2U]$ or more
independent trials of an experiment with success probability $\beta$, the probability of obtaining at least one success
exceeds $1-\frac{\varepsilon}{8}$. For $h,w\in \NN$ with $h,w>
100r$, define
$$
\Theta^R(h,w):=\{t:\lceil-r,r\rfloor\times 0
\text{~~is~joined~to~}\lceil w+h\ii,w\rfloor\times
t\text{~~within~}\lceil-w, w+h\ii\rfloor\}.
$$
And denote by $\mathbf{m}(\cdot)$ the Lebesgue measure on
$[0,\infty)$. Then $\mathbf{m}(\Theta^R(h,w))$ is the length of
infected time of the right side of the box $\lceil -w,
w+h\ii\rfloor$. Define $\Theta^L$, $\Theta^{UL}$, and $\Theta^{UR}$
similarly. Note that, for any $D\in\{L,R,UL,UR\}$ and $h,w\in\NN$,
\begin{equation}\label{e:5.33}
\{\Phi^D(h,w)=\emptyset\}=\{\Theta^D(h,w)=\emptyset\}.
\end{equation}
First, we will prove the following lemma.

\begin{lemma}\label{l:3.3}
One of the following two
assertions must be true.\\
$(1')$ There exist constants $h,w$ with $w=4h>100r$, such that
$$
\P(|\Phi^R(h,w)|+\mathbf{m}(\Theta^R(h,w))>U+V)>1-\frac{\varepsilon}{2},~~~\P(|\Phi^R(h,2w)|+\mathbf{m}(\Theta^R(h,2w))>U+V)>1-\frac{\varepsilon}{2}.
$$
$(2')$ There exist constants $h , w $ with $~ 8h\geq w$, such that
$$
\P(|\Phi^{UR}(h,w)|+\mathbf{m}(\Theta^{UR}(h,w))>U+V)>1-\frac{\varepsilon}{2},~~~
\P(|\Phi^R(2h,w)|+\mathbf{m}(\Theta^R(2h,w))>U+V)>1-\frac{\varepsilon}{2}.
$$
\end{lemma}

\pf Set $w_n=2^n$ for each $n>100r$. Since $\P(\xi^0$ dies out within
$\lceil -k, k+\infty\ii\rfloor)=1$ for all $k\in\NN$, we have, for every $n>100r$,
$$
\P(\xi^{\lceil-r,r\rfloor}\text{ dies out within } \lceil -w_n,
w_n+\infty\ii\rfloor)=1.
$$
This implies that we can find some $h_n\in
\{2^{w_n},~2^{w_n+1},~2^{w_n+2},\cdots\}$, such that
\begin{equation}\label{e:5.4}
\P(\xi^{\lceil-r,r\rfloor} \text{ dies out within } \lceil -w_n,
w_n+h_n\ii-\ii\rfloor)>1-\frac{\varepsilon^2}{8}.
\end{equation}
Without loss of generality, we suppose $(h_n)$ to be a strictly increasing
sequence. Then all sites being joined with $\lceil-w_n, w_n+h_n
\ii\rfloor$ are contained in $\lceil-w_{n+1}, w_{n+1}+h_{n+1}
\ii\rfloor$. By (\ref{e:5.4}), we have
\begin{equation}\label{e:5.5}\P(|\Phi^{UL}(h_n,w_n)|+|\Phi^{UR}(h_n,w_n)|=0~|~|\Phi(h_n,w_n)|+\mathbf{m}(\Theta(h_n,w_n))>2U+2V)>1-\frac{\varepsilon^2}{8}
\end{equation}
for  all $n>100r$.  For $h,w\in\NN$ with $h,w>100r$, denote
$$
\Theta(h,w):=\Theta^R(h,w)\cup\Theta^L(h,w)\cup\Theta^{UR}(h,w)\cup\Theta^{UL}(h,w).
$$
As before, let $\FF_n$ be the $\sigma$-field generated by the graphical representation within  $\lceil-w_n+h_n\ii,w_n+h_n\ii\rfloor~(n=1,2,\cdots)$. Note that, for any $n\in\NN$, if there is no flow passing through the edges $$\Xi(h_n,w_n):=\{(x,y)\in\EE:~x\in \Phi(h_n,w_n),
y\not\in\lceil-w_n,w_n+h_n\ii\rfloor\}$$ for every $t\in\Theta(h_n,w_n)$, then
$\xi^{\lceil-r,r\rfloor}$ must die out, since
 no  sites outside  $\lceil-w_n, w_n+h_n\ii\rfloor$ can be infected. Here, $\Phi(\cdot,\cdot)$ is defined as in (\ref{e:3.11}). Note that $|\Xi(h_n,w_n)|=|\Phi(h_n,w_n)|+2$ for any $n\in\NN$. And, for any $n\in\NN$, $A\subseteq\lceil-w_n,-w_n+h_n\ii\rfloor\cup\lceil-w_n+h_n\ii,w_n+h_n\ii\rfloor\cup\lceil w_n+h_n\ii,w_n\rfloor$, and $B\subseteq[0,\infty)$, we have $\Phi(h_n,w_n),\Theta(h_n,w_n)\in\FF_n$, and
\begin{align*}
&\P(\text{there is no flow passing through the edges in }\Xi(h_n,w_n),~\Phi(h_n,w_n)=A,~\Theta(h_n,w_n)=B~|~\FF_n)\\
\geq&\textbf{1}_{\{\Phi(h_n,w_n)=A,~\Theta(h_n,w_n)=B\}}\cdot[\E(\exp\{-\mathbf{m}(B)\cdot\xi\})]^{|A|+2},
\end{align*}
where $\xi$ is a random variable with law $\mu$. So
$$\P\left(\left.\xi^{\lceil-r,r\rfloor} \text{ dies
out}~\right|~\FF_n\right)\geq[\LLL(\mathbf{m}(\Theta(h_n,w_n)))]^{|\Phi(h_n,w_n)|+2}$$ for any $n\in\NN$, where $\LLL(t):=\E e^{-t\xi}$ is the Laplace transform of the random variable $\xi$. By the martingale convergence theorem,
$$\P\left(\left.\xi^{\lceil-r,r\rfloor} \text{ dies
out}~\right|~\FF_n\right)\rightarrow\textbf{1}_{\{\xi^{\lceil-r,r\rfloor} \text{ dies
out}\}}~~\text{a.s.}$$ as $n$ tends to infinity. So
$$\lim\limits_{n\rightarrow\infty}[\LLL(\mathbf{m}(\Theta(h_n,w_n)))]^{|\Phi(h_n,w_n)|+2}=0~\text{almost surely on}~\{\xi^{\lceil-r,r\rfloor}\text{ survives}\}.$$
But $\lim\limits_{n\rightarrow\infty}[\LLL(\mathbf{m}(\Theta(h_n,w_n)))]^{|\Phi(h_n,w_n)|+2}=0$ implies that $\lim\limits_{n\rightarrow\infty}[|\Phi(h_n,w_n)|+\mathbf{m}(\Theta(h_n,w_n))]=\infty$. So
$$\lim\limits_{n\rightarrow\infty}[|\Phi(h_n,w_n)|+\mathbf{m}(\Theta(h_n,w_n))]=\infty~\text{almost surely on}~\{\xi^{\lceil-r,r\rfloor}\text{ survives}\}.$$
Therefore,
$$
\P(\exists m,~\forall n>m,~|\Phi(h_n,w_n)|+\mathbf{m}(\Theta(h_n,w_n))>2U+2V~|~\xi^{\lceil-r,r\rfloor}\text{ survives} )=1.
$$
Hence there exists $n_0>100r$ such that, for $n>n_0$,
\begin{equation}\label{e:5.52}
\P(|\Phi(h_n,w_n)|+\mathbf{m}(\Theta(h_n,w_n))>2U+2V~|~\xi^{\lceil-r,r\rfloor}\text{ survives} )>1-\frac{\varepsilon^2}{16}.
\end{equation}
By (\ref{e:5.3}) and (\ref{e:5.52}), we get, for $n>n_0$,
\begin{equation}\label{e:5.53}
\P(|\Phi(h_n,w_n)|+\mathbf{m}(\Theta(h_n,w_n))>2U+2V)>1-\frac{\varepsilon^2}{8}.
\end{equation}
By (\ref{e:5.33}), (\ref{e:5.5}) and (\ref{e:5.53}), we have, for large $n$,
$$
\P(|\Phi^R(h_n,w_n)|+|\Phi^L(h_n,w_n)|+\mathbf{m}(\Theta^R(h_n,w_n))+\mathbf{m}(\Theta^L(h_n,w_n))>
2U+2V)>1-\frac{\varepsilon^2}{4}.
$$
Using the FKG inequality and the symmetric property again, we have
\begin{equation}\label{e:5.8}
\P(|\Phi^R(h_n,w_n)|+\mathbf{m}(\Theta^R(h_n,w_n))>
U+V)>1-\frac{\varepsilon}{2}
\end{equation}
for any sufficient large $n$. By (\ref{e:5.8}), we can conclude that one of the following two
assertions must be true.\\
$(1')$ There exist constants $r,h,w$ with $w=4h$, such that
$$
\P(|\Phi^R(h,w)|+\mathbf{m}(\Theta^R(h,w))>U+V)>1-\frac{\varepsilon}{2},~~~\P(|\Phi^R(h,2w)|+\mathbf{m}(\Theta^R(h,2w))>U+V)>1-\frac{\varepsilon}{2}.
$$
$(2')$ There exist constants $h , w $ with $~ 8h\geq w$, such that
$$
\P(|\Phi^{UR}(h,w)|+\mathbf{m}(\Theta^{UR}(h,w))>U+V)>1-\frac{\varepsilon}{2},~~~
\P(|\Phi^R(2h,w)|+\mathbf{m}(\Theta^R(2h,w))>U+V)>1-\frac{\varepsilon}{2}.
$$
The argument is a little modification from the proof of Case 1 to reduce the height, and is omitted here. We have finished the proof of the lemma.\qed\\

Comparing Lemma \ref{l:3.3} with Case 2, we
only need to prove the following.\\
(a) If $h$ and $w$ satisfy $(1')$, then
\begin{equation}\label{e:5.9}\P(|\Phi^{R}(h+N,w+4N)|>N~|~|\Phi^R(h,w)|+\mathbf{m}(\Theta^R(h,w))>U+V)>1-\frac{\varepsilon}{2}\end{equation}
and
\begin{equation}\label{e:5.91}\P(|\Phi^{R}(h+N,2w+8N)|>N~|~|\Phi^R(h,2w)|+\mathbf{m}(\Theta^R(h,2w))>U+V)>1-\frac{\varepsilon}{2}.\end{equation}
(b) If $h$ and $w$ satisfy $(2')$, then
\begin{equation}\label{e:5.92}\P(|\Phi^{UR}(h+N,w+8N)|>N~|~|\Phi^{UR}(h,w)|+\mathbf{m}(\Theta^{UR}(h,w))>U+V)>1-\frac{\varepsilon}{2}\end{equation}
and
\begin{equation}\label{e:5.93}\P(|\Phi^{UR}(h+N,2w+16N)|>N~|~|\Phi^{UR}(h,2w)|+\mathbf{m}(\Theta^{UR}(h,2w))>U+V)>1-\frac{\varepsilon}{2}.\end{equation}
We only prove (\ref{e:5.9}), since the proofs of (\ref{e:5.91})--(\ref{e:5.93}) are similar. Note that, if
\begin{equation}\label{e:5.10}
\P(|\Phi^{R}(h+N,w+4N)|>N,~|\Phi^{R}(h,w)|>U)\geq\left(1-\frac{\varepsilon}{4}\right)\cdot\P(|\Phi^{R}(h,w)|>U)
\end{equation}
and
\begin{align}\label{e:5.11}
&\P(|\Phi^{R}(h+N,w+4N)|>N,~|\Phi^{R}(h,w)|<U,~\mathbf{m}(\Theta^{R}(h,w))>V)\nonumber\\
\geq&\left(1-\frac{\varepsilon}{4}\right)\cdot\P(|\Phi^{R}(h,w)|<U,~\mathbf{m}(\Theta^{R}(h,w))>V),
\end{align}
then (\ref{e:5.9}) holds. Therefore, to prove (\ref{e:5.9}), it suffices to prove (\ref{e:5.10}) and (\ref{e:5.11}).\\

\noindent\textbf{Proof of (\ref{e:5.10})}\quad Let $h$ and $w$ satisfy $(1')$. Let $t_1$ be
the first time  that a site in $\lceil w, w+(h-2N)\ii\rfloor$ is
infected. That is,
$$
t_1:=\inf\{t:
 \lceil-r,r\rfloor\times 0\text{ is joined to }
\lceil w, w+(h-2N)\ii\rfloor\times t \text
{~within~}\lceil-w,w+h\ii\rfloor\times[0,\infty)\}.
$$
If $t_1<\infty$, then with probability 1, there exists a unique
infected site $x_1\in \lceil w, w+(h-2N)\ii\rfloor$ such that
$\lceil-r,r\rfloor\times 0\text{ is joined to } x_1\times t_1 \text
{~within~}\lceil-w,w+h\ii\rfloor\times[0,\infty)$. Generally, let
$t_k$ be the first time that a site in $\lceil w,
w+(h-2N)\ii\rfloor\backslash (\cup_{i=1}^{k-1} \lceil x_i-2N\ii,
x_i+N\ii\rfloor)$ is infected, and let $x_k$ be the corresponding
infected site if $t_k <\infty$. Denote by $E_k$ the event that
$x_k\times t_k$ is joined to every site of $\lceil x_k+4N,
x_k+4N+N\ii\rfloor\times (t_k+1)$ within $\{x_k\}\cup\lceil
x_k+1,x_k+4N+N\ii\rfloor$. If $E_k$ occurs, then
$|\Phi^R(h+N,w+4N)|>N$.
 By transitivity and rotation invariance of the space, we know that
 $(\textbf{1}_{E_k}|t_k<\infty)_{k=1}^{\infty}$ has the same distribution as $\textbf{1}_E$, where $E$ is defined in (\ref{e:5.31}).  Let
$$Y_k=\left\{\begin{array}{ll} \textbf{1}_{E_k},~~~~~~~~~~~~~~~~~~~~~~~~~~~~~~~~~~~~~~~~~~~~~~~~~~~~~~~~~~~~~~~~~~~~~~~~~~~~~~~~~~~~~\; \text{ if } t_k<\infty,\\
\text{an independent random variable with the same distribution as }\textbf{1}_E,
~~~~\text{ if } t_k=\infty.\end{array}\right.$$ Then $\P(Y_k=1)=1-\P(Y_k=0)=\alpha$.

Note that $Y_1,Y_2,\cdots$ are independent with respect to $\P$, since they are measurable with respect to the $\sigma$-fields generated by the graphical representations within mutually disjoint edge sets. Also, there exists $t_1<\cdots<t_{[U/20N]}<\infty$ almost surely if
$|\Phi^R(h,w)|>U$. Moreover, $\{|\Phi^R(h,w)|>U\}$ and
$\{\sum_{k=1}^{[U/20N]} Y_k\ge 1\}$ are increasing events.
Therefore, by the FKG inequality,
\begin{align*}
\P(|\Phi^{R}(h+N,w+4N)|>N,~|\Phi^R(h,w)|>U )&\geq\P(\text{some } E_k \text{ occurs },~|\Phi^R(h,w)|>U )\\
&=\P\left(|\Phi^R(h,w)|>U, \sum_{k=1}^{[U/20N]} Y_k\ge 1\right)\\
&\ge \P(|\Phi^R(h,w)|>U)\cdot \P\left( \sum_{k=1}^{[U/20N]}
Y_k\ge 1\right)\\
&\ge \left(1-\frac{\varepsilon}{4}\right)\cdot\P(|\Phi^R(h,w)|>U).
\end{align*}
Then (\ref{e:5.10})
holds, as desired.\qed\\

\noindent\textbf{Proof of (\ref{e:5.11})}\quad For any $x\in\lceil w,w+h\ii\rfloor$, set
$$T(x):=\mathbf{m}(\{t\geq0:~\lceil-r,r\rfloor\times0\text{~is~joined~to~}x\times t\text{~within~}\lceil-w,w+h\ii\rfloor\})$$ to
be the Lebesgue measure of the total infection time of $x$. So, if $|\Phi^R(h,w)|<U$ and
$\mathbf{m}(\Theta^R(h,w))>V$, then there exists $x\in\lceil
w,w+h\ii \rfloor$ such that $T(x)>\frac{V}{U}$. Define random events
\begin{align*}
A_0&=\left\{T(w)>\frac{V}{U}\right\},\\
A_1&=\left\{T(w)\leq\frac{V}{U},~T(w+\ii)>\frac{V}{U}\right\},\\
A_2&=\left\{T(w)\leq\frac{V}{U},~T(w+\ii)\leq\frac{V}{U},~T(w+2\ii)>\frac{V}{U}\right\},\\
&\cdots\cdots\cdots\cdots\cdots\cdots\cdots\cdots\\
A_h&=\left\{T(w)\leq\frac{V}{U},~\cdots,~T(w+(h-1)\ii)\leq\frac{V}{U},~T(w+h\ii)>\frac{V}{U}\right\}.
\end{align*}

For any $0\leq k\leq h$, suppose that $A_k$ occurs. We set $s_0=0$ and
$$s_i=\inf\{t\in(s_{i-1}+1,\infty):~\mathbf{m}(\{s_{i-1}+1<s<t:~x_k\text{~is~infected~at~time~}s\})=1\}$$
for $i=1,2,\cdots$ inductively. (Here, $\inf\emptyset$ is defined to be $+\infty$.) Define
$$D_i:=(s_{i-1}+1,s_i)\cap\{t\geq0:~x_k\text{~is~infected~at~time~}t\}$$
for $i=1,2,\cdots$. Then $\mathbf{m}(D_i)=1$ if $s_i<+\infty$. And
$d(D_{i-1},D_i)\geq1$ for $i=1,2,\cdots$, where
$$d(A,B):=\inf\{|a-b|:~a\in A,~b\in B\}$$ for any $A,B\subset\RR$.
Furthermore, for $i=1,2,\cdots$ define $$\tau_i:=\inf\{t\in
(s_{i-1}+1,s_i):~x_k\text{~is~infected~at~time~}t\}.$$ Note that $s_i<\infty$ implies that $\tau_i<\infty$ for $i=1,2,\cdots$.

 Denote by $F_i$ the event that
$x_k\times\tau_i$ is joined to every site of $\lceil x_k+4N,
x_k+4N+N\ii\rfloor\times(\tau_i+1)$ within $\{x_k\}\cup\lceil
x_k+1,x_k+4N+N\ii\rfloor$. If $F_i$ occurs, then
$|\Phi^R(h+N,w+4N)|>N$.
 By transitivity and rotation invariance of the space, we know that
 $(\textbf{1}_{F_i}|s_i<\infty)_{i=1}^{\infty}$ has the same distribution as $\textbf{1}_F$, where $F$ is the event that
$0\times0$ is joined to every site of $\lceil 4N,
4N+N\ii\rfloor\times1$ within $\{0\}\cup\lceil 1,4N+N\ii\rfloor$.
Let
$$Z_i=\left\{\begin{array}{ll} \textbf{1}_{F_i},~~~~~~~~~~~~~~~~~~~~~~~~~~~~~~~~~~~~~~~~~~~~~~~~~~~~~~~~~~~~~~~~~~~~~~~~~~~~~~~~~~~~~~\; \text{ if } s_i<\infty,\\
\text{an independent random variable with the same distribution as }\textbf{1}_F,
~~~~\text{ if } s_i=\infty.\end{array}\right.$$
By the strong Markov property under the quenched law, we know that $Z_1,Z_2,\cdots$ are independent with respect to $\P_{\lambda}$ for any fixed environment $\lambda$. And for any environment $\lambda$ such that $\lambda_e\geq a$ for all $e\in\{x_k\}\cup\lceil x_k+1,x_k+4N+N\ii\rfloor$, we have
$$\P_{\lambda}(Z_i=1)=1-\P_{\lambda}(Z_i=0)\geq\beta$$
by the monotonicity of the contact process. So, by our choice of $V$ and $U$,
\begin{equation}\label{e:5.12}
\P_{\lambda}\left(\sum\limits_{i=1}^{[V/2U]}Z_i\geq1\right)\geq1-\frac{\varepsilon}{8}.
\end{equation}
Turning to the annealed law, we get from (\ref{e:5.32}) and (\ref{e:5.12}) that
$$\P\left(\sum\limits_{i=1}^{[V/2U]}Z_i\geq1\right)\geq\P^{\mu}(\lambda_e\geq a\text{ for all }e\in\{x_k\}\cup\lceil x_k+1,x_k+4N+N\ii\rfloor)\cdot\left(1-\frac{\varepsilon}{8}\right)\geq1-\frac{\varepsilon}{4}.$$
Furthermore, note that there exists $s_1<\cdots<s_{[V/2U]}<\infty$ almost surely if
$|\Phi^R(h,w)|<U$ and $\mathbf{m}(\Theta^R(h,w))>V$.
 Therefore,
\begin{align*}
&\P(|\Phi^{R}(h+N,w+4N)|>N,~|\Phi^R(h,w)|<U,~\mathbf{m}(\Theta^R(h,w))>V,~A_k )\\
\geq&\P(\text{some } F_i \text{ occurs},~|\Phi^R(h,w)|<U,~\mathbf{m}(\Theta^R(h,w))>V,~A_k )\\
=&\P\left(|\Phi^R(h,w)|<U,~\mathbf{m}(\Theta^R(h,w))>V,~A_k,~\sum_{i=1}^{[V/2U]} Z_i\ge 1\right)\\
=&\P(|\Phi^R(h,w)|<U,~\mathbf{m}(\Theta^R(h,w))>V,~A_k)\cdot
\P\left( \sum_{i=1}^{[V/2U]}
Z_i\ge 1\right)\\
\geq&\left(1-\frac{\varepsilon}{4}\right)\cdot\P(|\Phi^R(h,w)|<U,~\mathbf{m}(\Theta^R(h,w))>V,~A_k).
\end{align*}
Here, the third equality holds because the event $\{|\Phi^R(h,w)|<U,~\mathbf{m}(\Theta^R(h,w))>V,~A_k\}$ is measurable with respect to the $\sigma$-field generated by the graphical representation within $G_1=\lceil-w,w+h\ii\rfloor$, while the event $\left\{ \sum\limits_{i=1}^{[V/2U]}
Z_i\ge 1\right\}$ is measurable with respect to the $\sigma$-field generated by the graphical representation within $G_2=\{x_k\}\cup\lceil x_k+1,x_k+4N+N\ii\rfloor$. The two events are independent, since $G_1$ and $G_2$ are disjoint edge sets which share no common edges. Next, note that $$\{|\Phi^R(h,w)|<U\}\cap\{\mathbf{m}(\Theta^R(h,w))>V\}\subseteq\bigcup\limits_{k=0}^{h}A_k,$$
and $A_k,~k=0,1,\cdots,h$, are mutually exclusive events.
Therefore,
\begin{align*}
&\P(|\Phi^{R}(h+N,w+4N)|>N,~|\Phi^R(h,w)|<U,~\mathbf{m}(\Theta^R(h,w))>V)\\
=&\sum\limits_{k=0}^h\P(|\Phi^{R}(h+N,w+4N)|>N,~|\Phi^R(h,w)|<U,~\mathbf{m}(\Theta^R(h,w))>V,~A_k )\\
\geq&\left(1-\frac{\varepsilon}{4}\right)\cdot\sum\limits_{k=0}^h\P(|\Phi^R(h,w)|<U,~\mathbf{m}(\Theta^R(h,w))>V,~A_k)\\
=&\left(1-\frac{\varepsilon}{4}\right)\cdot\P(|\Phi^R(h,w)|<U,~\mathbf{m}(\Theta^R(h,w))>V)
\end{align*}
Then (\ref{e:5.11}) holds, as desired.\qed\\

All the above arguments together lead to the proof of Case 2.

\subsection{Proof of Case 3}

In this subsection, we shall prove that the block conditions hold if $\xi^0$ survives within some \lq slab\rq~ with positive probability. Choose fixed $K\in \NN$ such that
\begin{equation}\label{e:6.2}\P\left(\xi^0 \text{~survives~within~}  \lceil
-K, K+\infty\ii\rfloor\right)=c>0.\end{equation}
For any $x\in \HH$, $m,n\in\NN$, and $t>0$, denote by $A(x,t,
m,n)$ the event that  $x\times t$ is joined to $\lceil
x+m+n\ii+2K+2K\ii,x+m+n\ii+4K+4K\ii\rfloor\times[t,\infty)$ within
$\{x\}\cup\lceil x+\ii+K, x-K+2K\ii+n\ii\rfloor\cup\lceil
x-K+2K\ii+n\ii, x+m+n\ii+4K+4K\ii\rfloor$. See Figure 2 for intuition.
Then, for  $m\in \NN$, $t>0$, and  $x\in \HH$ with $\Im(x)>K$, define
$$ T(x,t,m):=\inf\{s\geq t:~x\times t \text{~~is~ joined~
to~}\lceil x-K\ii+m, x+K\ii+m\rfloor\times s
\text{~within~}\{x\}\cup\lceil x+1-K\ii, x+K\ii+m\rfloor\}.
$$
And similarly, define
$$
T(x,t,m\ii):=\inf\{s\geq t:~x\times t \text{~~is~ joined~ to~}\lceil
x-K+m\ii, x+K+m\ii\rfloor\times s \text{~within~}\{x\}\cup\lceil
x+\ii-K, x+K+m\ii\rfloor\}.
$$
for $x\in\HH$, $m\in\NN$ and $t>0$. See Figure 3 for intuition.
 \begin{figure}[H]\label{picc:4}
 \center
 \includegraphics[width=12.0true cm, height=9.6true cm]{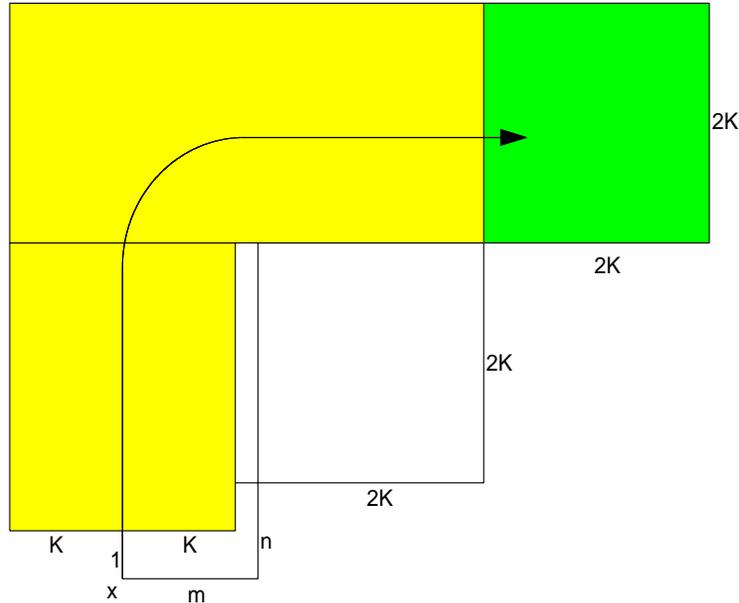}
 \caption{Description of $A(x,t,m,n)$}
 \end{figure}

 \begin{figure}[H]\label{picc:5}
 \center
 \includegraphics[width=12.0true cm, height=9.6true cm]{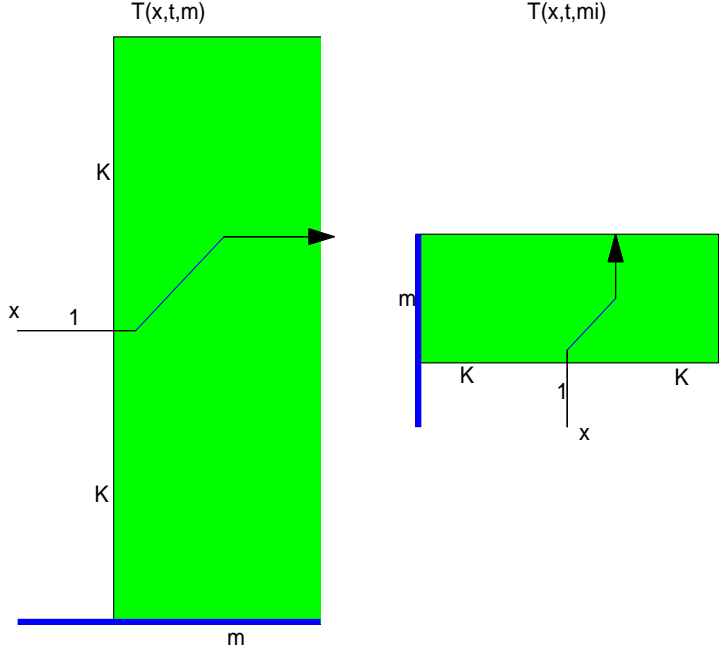}
 \caption{Description of $T(x,t,m)$ and $T(x,t,m\ii)$}
 \end{figure}
We then have the following lemma, which is essential to the proof of Case 3.
\begin{lemma}\label{l:3.4}
There exists $\alpha>0$ which is independent of $x,t,m$ and $n$, such that
\begin{equation}\label{e:6.3} \P(A(x,t,
m,n))>\alpha.
\end{equation}
\end{lemma}
\pf For $x\in\HH$ and $t>0$, denote by $C(x,t)$ the event that $x\times t$ is joined  to $(x+3K+3K\ii)\times
(t+1)$ within $\lceil x,x+3K\ii\rfloor\cup\lceil x+3K\ii,
x+3K+3K\ii\rfloor$. By
translation invariance we have that $\P(C(x,t))=\P(C(0,0))$ for any $x\in\HH$ and $t>0$. We next prove that $$\alpha:=\frac{c^2}{4}\cdot\P(C(0,0))$$ satisfies (\ref{e:6.3}), where $c$ is the positive constant as defined in (\ref{e:6.2}). By (\ref{e:6.2}) and the translation invariance, we have, for any $x\in\HH$, $m\in\NN$, and $t>0$,
\begin{equation}\label{e:6.6}
\P(T(x,t,m\ii)<\infty)=\P(T(0,0,h\ii)<\infty)\geq\P\left(~\xi^0 \text{~survives~within~}  \lceil
-K, K+\infty\ii\rfloor\right)=c.
\end{equation}
Furthermore, by rotation invariance, for any $m\in\NN$, $t>0$, and $x\in\HH$ with $\Im(x)>K$,
\begin{equation}\label{e:6.7}
\P(T(x,t,m)<\infty)=\P(T(x,t,m\ii)<\infty)\geq c.
\end{equation}
Next, if
$T(x,t,n\ii)<\infty$, then let $X(x,t,n\ii)$ be the corresponding
infected site.  For $x\in \HH$, $m,n\in\NN$, and $t>0$, define
$$
D(x,t,m,n):=\{T(x,t,n\ii)<\infty\}\cap
C(X(x,t,n\ii),T(x,t,n\ii))\cap\{T(X(x,t,n\ii)+3K+3K\ii,
T(x,t,n\ii)+1,m)<\infty\}.
$$
Obviously, for any $x\in\HH$, $t>0$, and $m,n\in\NN$, we have
\begin{equation}\label{e:6.8} D(x,t, m,n)\subseteq
A(x,t,m,n).\end{equation}
Let $\FF$ denote the $\sigma$-field generated by the graphical representation within $\{x\}\cup\lceil x-K+\ii,x+K+n\ii\rfloor$. Then $X(x,t,n\ii)$ and $T(x,t,n\ii)$ are measurable with respect to $\FF$. So, by (\ref{e:6.6})--(\ref{e:6.8}), we have
\begin{align*}
&\P(A(x,t,m,n))\geq\P( D(x,t,m,n) )\\
=&\E(\P(T(x,t,n\ii)<\infty,~C(X(x,t,n\ii),T(x,t,n\ii)),~T(X(x,t,n\ii)+3K+3K\ii,T(x,t,n\ii)+1,m)<\infty~|~\FF))\\
=&\E(\P(s<\infty,~C(y,s),~T(y+3K+3K\ii,s+1,m)<\infty~|~\FF)~|_{y=X(x,t,n\ii),~s=T(x,t,n\ii)})\\
=&\E(\mathbf{1}_{\{s<\infty\}}\cdot\P(C(y,s))\cdot\P(T(y+3K+3K\ii,s+1,m)<\infty)~|_{y=X(x,t,n\ii),~s=T(x,t,n\ii)})\\
=&\P(T(x,t,n\ii)<\infty)\cdot\P(C(0,0))\cdot\P(T(0,0,m)<\infty)\\
\geq&c^2\cdot\P(C(0,0))>\alpha
\end{align*}
We next explain the third equality in detail. By definition, $X(x,t,n\ii)$ takes a value in $\lceil x-K+n\ii,x+K+n\ii\rfloor$. For any fixed $y\in\lceil x-K+n\ii,x+K+n\ii\rfloor$ and $s>0$, the event $C(y,s)$ is measurable with respect to the $\sigma$-field generated by the graphical representation within $G_1=\lceil y,y+3K\ii\rfloor\cup\lceil y+3K\ii,y+3K+3K\ii\rfloor$, while the event $\{T(y+3K+3K\ii,s+1,m)<\infty\}$ is measurable with respect to the $\sigma$-field generated by the graphical representation within $G_2=\{y+3K+3K\ii\}\cup\lceil y+3K+1+2K\ii,y+3K+m+4K\ii\rfloor$. Note that $G_1$ and $G_2$ are disjoint with $\{x\}\cup\lceil x-K+\ii,x+K+n\ii\rfloor$, respectively. As a result, the events $C(y,s)$ and $\{T(y+3K+3K\ii,s+1,m)<\infty\}$ are independent of $\FF$, respectively. Furthermore, the two events are independent since $G_1$ and $G_2$ are disjoint edge sets which share no common edges.

From the above arguments, we get the inequality in (\ref{e:6.3}). Therefore, we have completed the proof of Lemma \ref{l:3.4}.\qed\\

\noindent\textbf{Proof of Case 3}\quad Fix $N\in\NN$ and $\varepsilon>0$ sufficiently small. Let $N_1$ be large
enough to ensure that, in $N_1$ or more independent trials of an experiment with
success probability $\alpha$, the probability of obtaining at least
$N$ success exceeds $1-\varepsilon/2$. Here, $\alpha$ is the positive constant in Lemma \ref{l:3.4}. By Lemma \ref{l:3.1}, there
exists $r$ such that $\xi^{\lceil -r, r\rfloor}$ survives with
probability greater than $1-\varepsilon/4$. For $m,n\in\NN$ and $t>0$, define
\begin{align*}
&B(m,n,t,UR);=\{\lceil -r, r\rfloor\times0\text{ is joined to }z\times t\text{ for all }z\in \lceil n+m\ii, n+m\ii+3KN_1\rfloor\text{ within }B_0(m)\},\\
&B(m,n,t,R):=\{\lceil -r,r\rfloor\times0\text{ is joined to }z\times t\text{ for all }z\in \lceil m+n\ii, m+n\ii+3KN_1\ii\rfloor\text{ within }B_0(m)\},\\
&B(m,n,t,UL);=\{\lceil -r, r\rfloor\times0\text{ is joined to }z\times t\text{ for all }z\in \lceil -n+m\ii, -n+m\ii-3KN_1\rfloor\text{ within }B_0(m)\},\\
&B(m,n,t,L):=\{\lceil -r,r\rfloor\times0\text{ is joined to }z\times t\text{ for all }z\in \lceil -m+n\ii, -m+n\ii+3KN_1\ii\rfloor\text{ within }B_0(m)\}.
\end{align*}
See Figure 4 for intuition. Then,
define
$$\rho:=\inf\{m\in\NN:~\exists t>0,~n\in\NN~\text{and}~D\in\{UR,R,UL,L\}~\text{s.t.}~B(m,n,t,D)\text{~occurs}\}.$$

 \begin{figure}[H]\label{picc:6}
 \center
 \includegraphics[width=12.0true cm, height=9.0true cm]{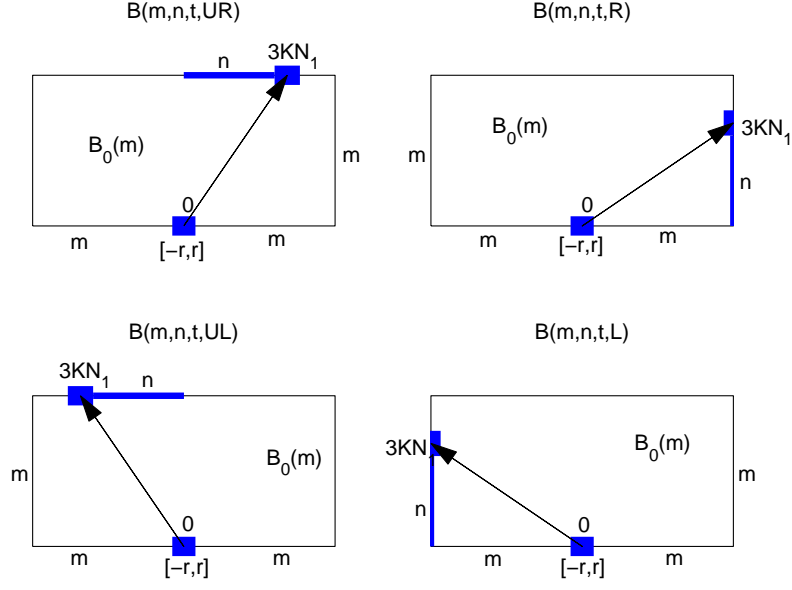}
 \caption{Description of $B(m,n,t,UR)$, $B(m,n,t,R)$, $B(m,n,t,UL)$ and $B(m,n,t,L)$}
 \end{figure}

We next prove that
\begin{equation}\label{e:6.81}
\P(\rho<\infty~|~\xi^{\lceil-r,r\rfloor}\text{ survives})=1.
\end{equation}
Define
$$p:=\P(\exists0<t<\infty,~\text{s.t.}~0\times0\text{ is joined to }z\times t\text{ for all }z\in\lceil0,3KN\ii\rfloor\text{ within }\lceil0,3KN\ii\rfloor).$$
Then $p>0$. For any $m\in\NN$, denote by $Z_m$ the first infected site in
$$\partial B_0(m):=\lceil-m,-m+m\ii\rfloor\cup\lceil-m+m\ii,m+m\ii\rfloor\cup\lceil m+m\ii,m\rfloor,$$
and by $\tau_m$ the corresponding infected time. Furthermore, denote by $\GG_m$ the $\sigma$-field generated by the graphical representation within $\lceil-m,m+m\ii\rfloor\setminus\partial B_0(m)$. If $\xi^{\lceil-r,r\rfloor}$ survives, then $\tau_m<\infty$ and $\tau_m,Z_m\in\GG_m$ for any $m\in\NN$. Next, for any $m\geq8KN_1$, we divide $\partial B_0(m)$ into four parts as follows:
$$\partial B_0^L(m):=\lceil-m,-m+m\ii\rfloor,~~~~~~~~\partial B_0^{UL}(m):=\lceil-m+m\ii,m\ii\rfloor,$$
$$\partial B_0^{UR}(m):=\lceil m\ii,m+m\ii\rfloor,~~~~~~~~\partial B_0^{R}(m):=\lceil m+m\ii,m\rfloor.$$
Since $m\geq8KN_1$, no matter which part of $D\in\{L,UL,UR,R\}$ that $Z_m$ lies in, we can find a seed $S_m$ with length $3KN_1$ such that $S_m$ lies entirely in the same part as $Z_m$ and one endpoint of $S_m$ is $Z_m$. For $m\geq8KN_1$, denote
$$A_m:=\{\exists\tau_m<t<\infty,~\text{s.t.}~Z_m\times\tau_m\text{ is joined to }z\times t\text{ for all }z\in S_m\text{ within }S_m\}.$$
Then, by translation and rotation invariance, if $\xi^{\lceil-r,r\rfloor}$ survives, then
$$\P(\rho<\infty~|~\GG_m)\geq\P(A_m~|~\GG_m)=p>0$$ for any $m\geq8KN_1$. That is,
\begin{equation}\label{e:6.82}
\P(\P(\rho<\infty~|~\GG_m)\geq p>0\text{ for all }m\geq8KN_1~|~\xi^{\lceil-r,r\rfloor}\text{ survives})=1.
\end{equation}
Furthermore, using the martingale convergence theorem, we can get that
$$\P(\rho<\infty~|~\GG_n)\rightarrow\textbf{1}_{\{\rho<\infty\}}~~\text{a.s.}$$ as $n$ tends to infinity. So, by (\ref{e:6.82}), we get
$$\P(\rho<\infty~|~\xi^{\lceil-r,r\rfloor}\text{ survives})=1.$$ Therefore, (\ref{e:6.81}) holds.

From (\ref{e:6.81}), we can get that there exists a positive integer
$N_2>100r$ such that
\begin{equation}\label{e:6.4}
\P(\rho<N_2~|~\xi^{\lceil-r,r\rfloor}\text{ survives})>1-\frac{\varepsilon}{4}.
\end{equation}
Set
$$
h:=4KN_1+N_2~~~~\text{and~~~~}w:=4h.
$$
We next prove that this choice of $r,h$, and $w$ satisfies (\ref{e:3.2}). If
$\rho<N_2$, let
$$
\tau:=\inf\{t>0:~\exists n\in\NN\text{ and }D\in\{L,UL,UR,R\}\text{ s.t. }B(\rho,n,t,D)
\text{~occurs}\}.
$$
Obviously, if $\rho<N_2$, then $\tau<\infty$. Let
$$
\gamma:=\inf\{n\in\NN:~\exists D\in\{L,UL,UR,R\}\text{ s.t. }B(\rho,n,\tau,D)\text{~occurs}\}.
$$
We divide our problem into four cases. (I):$B(\rho,\gamma,\tau,UR)$
occurs. (II):$B(\rho,\gamma,\tau,R)$ occurs. (III):$B(\rho,\gamma,\tau,LR)$
occurs. (IV):$B(\rho,\gamma,\tau,L)$ occurs. We only prove Cases (I) and (II), since Cases (III) and (IV) can be easily achieved by symmetry.

Case (I). Suppose that $B(\rho,\gamma,\tau,UR)$ occurs. Then $\lceil -r,
r\rfloor\times0$ is joined to $z\times t$  for all $z\in \lceil
\gamma+\rho\ii, \gamma+\rho\ii+3KN_1\rfloor$ within $B_0(\rho)$. For
$0\le j\le N_1-1$, let
$$
x_j:=\gamma+\rho\ii+K+3Kj \text{~~~and~~~} n_j=4KN_1-3Kj.
$$
Then $$\{x_j\}\cup\lceil x_j+\ii+K,
x_j-K+n_j\ii+2K\ii\rfloor\cup\lceil x_j-K+n_j\ii+2K\ii,
x+w+n_j\ii+4K+4K\ii\rfloor,~~0\le j\le N_1-1$$ are disjoint. By the
assumption of (\ref{e:6.3}) and the definition of $N_1$, with
probability greater than $1-\frac{\varepsilon}{2}$ there are at
least $N$ events in $\{A(x_j, \tau, w,n_j),~~0\le j\le N_1-1\}$ occur.
We can see that, if $A(x_j, \tau, w,n_j)$ occurs, then $\lceil -r,
r\rfloor\times 0$ is joined to $\lceil
w+\Im(x_j)\ii+n_j\ii+2K\ii,w+\Im(x_j)\ii+n_j\ii+4K\ii\rfloor\times[0,\infty)$
%
within
$$
B_0(\rho)\cup\lceil x_j+\ii+K, x_j-K+n_j\ii+2K\ii\rfloor\cup\lceil
x_j-K+n_j\ii+2K\ii, w+\Im(x_j)\ii+n_j\ii+4K\ii\rfloor\subset\lceil
-w, w+h\ii\rfloor.
$$
Therefore, conditioned on $B(\rho,\gamma,\tau,UR)$
occurs, the probability of $ |\Phi^R(h ,w)|>N$ is greater than
$1-\frac{\varepsilon}{2}$.

Case (II). Suppose that $B(\rho,\gamma,\tau,R)$ occurs. Then $\lceil -r,
r\rfloor\times0$ is joined to $z\times t$  for all $z\in \lceil
\gamma+\rho\ii, \gamma+\rho\ii+3KN_1\ii\rfloor$ within $B_0(\rho)$.
For $0\le j\le N_1-1$, let
$$
x_j:=\gamma+\rho\ii+K\ii+3Kj\ii.
$$
Then $$\{x_j\}\cup\lceil x_j+1-K\ii, x_j+w+K\ii\rfloor,~~0\le j\le
N_1-1$$ are disjoint. By the assumption of (\ref{e:6.3}) and the
definition of $N_1$, with probability greater than
$1-\frac{\varepsilon}{2}$ there are at least $N$ events in $\{A(x_j,
\tau, w,n_j),~0\le j\le N_1-1\}$ occur. We can see that if $A(x_j,
\tau, w,n_j)$ occurs, then $\lceil -r, r\rfloor\times 0$ is joined
to $\lceil
w+\Im(x_j)\ii-K\ii,w+\Im(x_j)\ii+K\ii\rfloor\times[0,\infty)$
%
within
$$
B_0(\rho)\cup\lceil x_j+1-K\ii, x_j+w+K\ii\rfloor\subset\lceil -w,
w+h\ii\rfloor.
$$
 Therefore, conditioned on $B(\rho,\gamma,\tau,R)$
occurs, the probability of $ |\Phi^R(h ,w)|>N$ is greater than
$1-\frac{\varepsilon}{2}$.

By the above analysis, we have that, conditioned
on $\rho<N_2$, the probability of  $ |\Phi^R(h ,w)|>N $ is greater
than $1-\frac{\varepsilon}{2}$. Together with (\ref{e:6.4}), we get
\begin{align*}
&\P(|\Phi^R(h,w)|>N)\\
\geq&\P(|\Phi^R(h,w)|>N~|~\rho<N_2)\cdot\P(\rho<N_2~|~\xi^{\lceil-r,r\rfloor}\text{ survives})\cdot\P(\xi^{\lceil-r,r\rfloor}\text{ survives})\\
>&\left(1-\frac{\varepsilon}{2}\right)\left(1-\frac{\varepsilon}{4}\right)\left(1-\frac{\varepsilon}{4}\right)>1-\varepsilon.
\end{align*}
Similarly, we can prove that $\P(|\Phi^R(h,2w)|>N)>1-\varepsilon$. So we have proved Case 3.\qed\\

The three subsections above give the whole proof of the \lq block conditions\rq~, Proposition \ref{p:3.2}. Next, we make further analysis. Let $G$ be the
event that $0\times 0$ is joined to every site of
$\lceil-r+4r\ii,r+4r\ii\rfloor\times 1$ within $\ang{-r,r+4r\ii}$.
Fix $N\ge \frac{20r\log\varepsilon}{\log(1-\P(G))}+1$
which is large enough to ensure that, in $[N/20r]$ or more
independent trials of an experiment with success probability
$\P(G)$, the probability of obtaining at least one success
exceeds
$1-\varepsilon$. We then have the next lemma.
\begin{lemma}\label{l:3.2} Suppose that $\P(|\Phi^R(h,w)|>N)>1-\varepsilon$.
Then, with $\P$-probability greater than $1-2\varepsilon$,
there exist $x\in \lceil w+4r, w+4r+h\ii\rfloor$ and $t>0$, such
that the horizontal seed $(0\times 0)_r$ is joined to the vertical
seed $(x\times t)_r$  within $\ang{-w-1,
w+4r+h\ii}.$\\
\end{lemma}
\pf Let $t_1$ be the first time  that some site in $\lceil w+2r\ii,
w+(h-2r)\ii\rfloor$ is infected. That is,
$$
t_1:=\inf\{t:
 \lceil-r,r\rfloor\times 0\text{ is joined to }
\lceil w+2r\ii, w+(h-2r)\ii\rfloor\times t \text
{~within~}\lceil-w,w+h\ii\rfloor\}.
$$
If $t_1<\infty$, then with probability 1 there exists a unique
infected site $x_1\in \lceil w+2r\ii, w+(h-2r)\ii\rfloor$ such that
$\lceil-r,r\rfloor\times 0\text{ is joined to } x_1\times t_1 \text
{~within~}\lceil-w,w+h\ii\rfloor$. Generally, let $t_k$ be the first
time that some site in $\lceil w+2r\ii, w+(h-2r)\ii\rfloor\backslash
(\cup_{i=1}^{k-1} \lceil x_i-3r\ii, x_i+3r\ii\rfloor)$ is infected,
and let $x_k$ be the corresponding infected site if $t_k <\infty$.
Denote by $G_k$ the event that $x_k\times t_k$ is joined to every
site of $\lceil x_k+4r-r\ii, x_k+4r+r\ii\rfloor\times (t_k+1)$
within $\ang{x_k-r\ii,x_k+4r+r\ii}$. If $G_k$ occurs, then the
horizontal seed $(0\times 0)_r$ is joined to the vertical seed
$(x_k\times t_k)_r$ within $\ang{-w-1, w+4r+h\ii}$.
By transitivity and rotation invariance of the space, we know that
$(\textbf{1}_{G_k}|t_k<\infty)$ has the same distribution as $\textbf{1}_G$.  Let
$$Y_k=\left\{\begin{array}{ll} \textbf{1}_{G_k},~~~~~~~~~~~~~~~~~~~~~~~~~~~~~~~~~~~~~~~~~~~~~~~~~~~~~~~~~~~~~~~~~~~~~~~~~~~~~~~~~~~~~\; \text{ if } t_k<\infty,\\
\text{an independent random variable with the same distribution as }\textbf{1}_G,
~~~~\text{ if } t_k=\infty.\end{array}\right.$$ Then
$\P(Y_k=1)=1-\P(Y_k=0)=\P(G)$.

Note that $Y_1,Y_2,\cdots$ are independent with respect to $\P$, since they are measurable with respect to the $\sigma$-fields generated by the graphical representations within mutually disjoint edge sets. Also, there exists $t_1<\cdots<t_{[N/20r]}<\infty$ almost
surely if $|\Phi^R(h,w)|>N$. Therefore,
\begin{align*}
\P(\text{some } G_k \text{ occurs} )&\ge
\P\left(|\Phi^R(h,w)|>N,~\sum_{k=1}^{[N/20r]} Y_k\ge 1\right)\\
&\ge \P(|\Phi^R(h,w)|>N)+ \P\left(
\sum_{k=1}^{[N/20r]}
Y_k\ge 1\right)-1\\
&\ge 1-2\varepsilon.
\end{align*}
So there exist $x\in \lceil w+4r+2r\ii, w+4r+(h-2r)\ii\rfloor$ and
$t>0$, such that the horizontal seed $(0\times 0)_r$ is joined to
the vertical seed $(x\times t)_r$ within $\ang{-w-1, w+4r+h\ii}$
with $\P$-probability greater than $1-2\varepsilon$.  \qed
\\

\noindent\textbf{Remark.} A similar conclusion holds for $\Phi^L$, $\Phi^{UR}$, and $\Phi^{UL}$.\\

Now, we give the following proposition, which is essential to the analysis in the following sections. See Figure 5 for intuition.
\begin{prop}\label{p:3.1}
Suppose that $\P(\xi^0~\text{survives})>0$. Then, for any $\varepsilon>0$ sufficiently small, there exist
$r\ge 1$ and $h\ge 100r$ such that the following three assertions hold with $\P$-probability greater than $1-\varepsilon$.\\
(i) The horizontal seed $(0\times0)_r$ is joined to a vertical seed
$(x\times t)_r$ within $\ang{ -4h-1,w+h\ii}$ for some $ 4h+4r\le w <
4.0001h,~ \Re(x)=w$, and
$t>0$.\\
(ii) The horizontal seed $(0\times0)_r$ is joined to a vertical seed
$(x\times t)_r$ within $\ang{ -8h-1,w+h\ii}$ for some $ 8h+4r\le
w<8.0001h,~ \Re(x)=w$, and $t>0$. \\
(iii) The horizontal seed $(0\times0)_r$ is joined to a vertical
seed $(x_1\times t_1)_r$  within $\ang{ -8h-1,w_1+h\ii}$ for some $
8h+4r\le w_1<8.0001h$ and $t_1>0$; and the horizontal seed
$(0\times0)_r$ is joined to a vertical seed $(x_2\times t_2)_r$
within $\ang{-w_2+h\ii,8h+1}$ for some $ 8h+4r\le w_2<8.0001h$ and
$t_2>0$.\\
\end{prop}

\begin{figure}[H]\label{picc:2}
 \center
 \includegraphics[width=14.0true cm, height=9.5true cm]{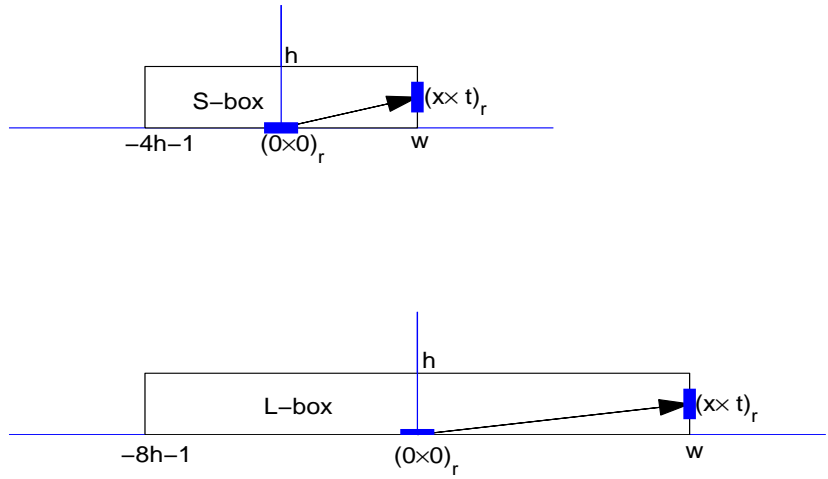}
 \caption{Construction of blocks}
 \end{figure}

\pf When $\P(\xi^0~\text{survives})>0$, either $(\verb"1")$ or $(\verb"2")$ of
Proposition \ref{p:3.2} is true. If $(\verb"1")$ is true, then, by
Lemma \ref{l:3.2},  $(i)$ and $(ii)$ hold. If $(\verb"2")$ is true,
we can prove the first two conclusions by iterating Lemma
\ref{l:3.2}; see Figure 6. Furthermore, by $(ii)$ together with the symmetric
property and the FKG inequality, we can get $(iii)$ in both cases.
So we have completed the proof of the proposition. \qed
\begin{figure}[H]\label{picc:3}
 \center
 \includegraphics[width=18.0true cm, height=2.9true cm]{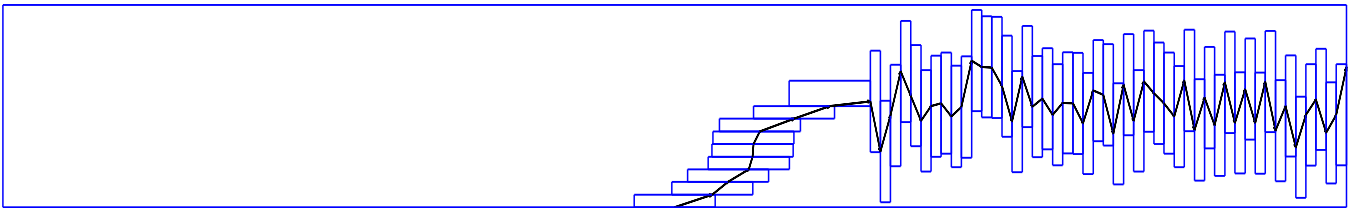}
 \caption{Construction of $(1)$ through  $(2)$}
 \end{figure}

\section{Dynamic renormalization}
From now on, for simplicity, we call the two kinds of edge sets displayed in Figure 5 \emph{S-boxes} and \emph{L-boxes}, respectively. (\lq S\rq~ stands for \lq short\rq; \lq L\rq~ stands for \lq long\rq.) Rigorously, S-boxes are edge sets having the same shape as $\ang{ -4h-1,w+h\ii}~(4h+4r\le w <
4.0001h,~ \Re(x)=w)$ described in Part (i) of Proposition \ref{p:3.1}, while L-boxes are edge sets having the same shape as $\ang{ -8h-1,w+h\ii}~(8h+4r\le
w<8.0001h,~ \Re(x)=w)$ described in Part (ii) of Proposition \ref{p:3.1}. The ratio of the width to the height in an S-box is nearly $8:1$, while the ratio of the width to the height in an L-box is nearly $16:1$. Translations and rotations are allowed. These edge sets are called \lq boxes\rq~ since the endpoints of each edge box form a rectangle on $\HH$. From Proposition \ref{p:3.1}, we are able to find some S-boxes and L-boxes such that, with large $\P$-probability, a horizontal seed on the bottom of each box is joined within the box to a vertical seed on the right. Figure 5 gives an intuition for it.\\

Next, we use these S-boxes and L-boxes to construct a route so that, with large probability, a seed in a fixed square is joined through the route to some seeds in the other two fixed squares (one above, the other on the right). The rigorous arguments are as follows. Set $\MM=10^7$ from now on. For any $\varepsilon>0$ sufficiently small, fix $r=r(\varepsilon)$ and $h=h(\varepsilon)$ satisfying Proposition \ref{p:3.1} henceforth. Next, for $x\in\HH$, $m\in\ZZ$, and $n\in\ZZ^+$, define
$$R_{m,n}(x):=\lceil
a+m\textrm{{M}}h+n\textrm{{M}}h\ii,
b+m\textrm{{M}}h+n\textrm{{M}}h\ii\rfloor=\lceil a, b\rfloor
+\textrm{{M}}h(m+n\ii),$$ where
$a=100h[\Re(x)/100h]+100h[\Im(x)/100h]\ii$ and $b=a+100(1+\ii)$.
Then $R_{m,n}(x)$ is a square and $x\in R_{0,0}(x)$.

Suppose that $(x\times s)_r$ is a seed (no matter whether it is horizontal or vertical). We next construct a
route by which this seed is joined to two vertical seeds in
$R_{0,1}(x)$ with large probability in the following way (see Figure 7 for intuition). Use S-boxes (horizontal and vertical
boxes alternatively) to let the seed spread in the northwest
(\lq~$\nwarrow$\rq) direction. If the infection surpasses the line
$\{y:~\Re(y)=\Re(a)+30h\}$, then use two L-boxes to change the
spread into the northeast (\lq$\nearrow$~\rq) direction. If the
infection surpasses the line $\{y:~\Re(y)=\Re(a)+70h\}$, then use
two L-boxes to change the spread into the northwest direction.
Iterate the procedure until the infection reaches $R_{0,1}(x)$. Then use an extra
L-box to get the two infected seeds we want. As a result, by the
route described above, the initial vertical seed $(x\times s)_r$ may be
joined to two vertical seeds $(y_1\times t_1)_r$ and $(y_2\times
t_2)_r$, where $y_1,y_2\in R_{0,1}(x)$. The vertical seed $(y_1\times t_1)_r$~(centering at $y_1$ and being generated at time $t_1$) will be used to make the next route in the \lq above\rq~ direction, while the vertical seed $(y_2\times t_2)_r$~(centering at $y_2$ and being generated at time $t_2$) will be used to make the next route in the \lq right\rq~ direction. See Figure 7 for the
precise positions of $y_1$ and $y_2$. Note that the route lies entirely in
$\lceil a, b+\MM h\rfloor$.

 \begin{figure}[H]\label{picc:7}
 \center
 \includegraphics[width=8.0true cm, height=15.5true cm]{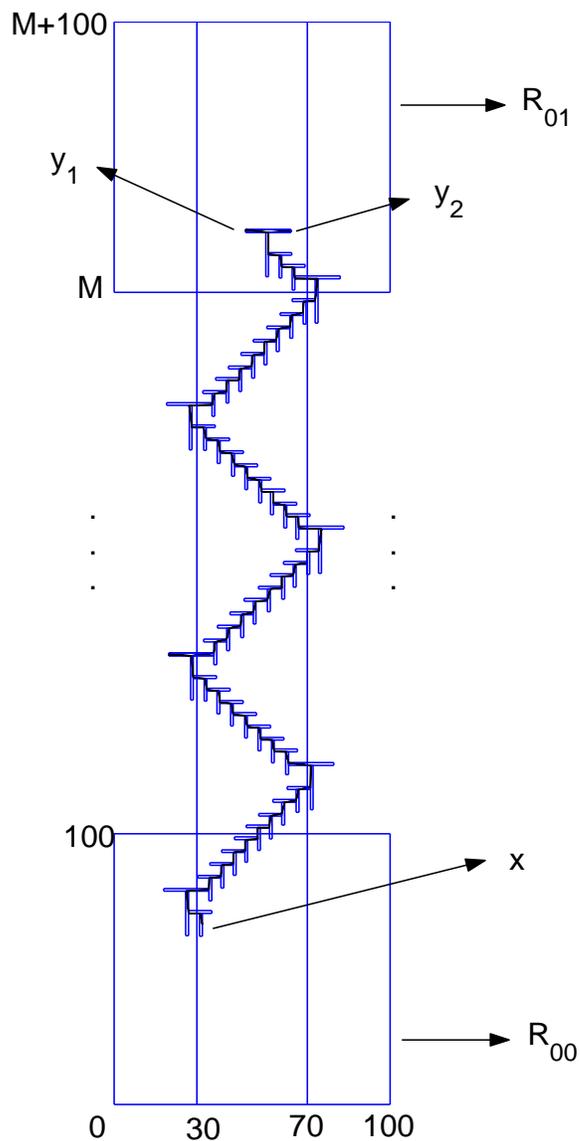}
 \caption{Producing new seeds in $R_{0,1}(x)$}
 \end{figure}

The number of steps in the above procedure is no more than $\MM$. So, by Proposition \ref{p:3.1} together with the fact that the events are independent if they are measurable with respect to $\sigma$-fields generated by graphical representations within disjoint subgraphs (this has been used several times in Section 3; for details readers can refer to Lemmas 3.5 and 3.6 of \cite{Chen-Yao2009}), we can get $t_1+t_2<\infty$ with
large probability. If $t_1+t_2<\infty$, then the above procedure generates two seeds as required. Similarly, we can construct a
route by which the seed $(x\times s)_r$ is joined to two horizontal seeds in
$R_{1,0}(x)$ with large probability. See Figure 8 for intuition.\\

  \begin{figure}[H]\label{picc:8}
 \center
 \includegraphics[width=14.0true cm, height=3.7true cm]{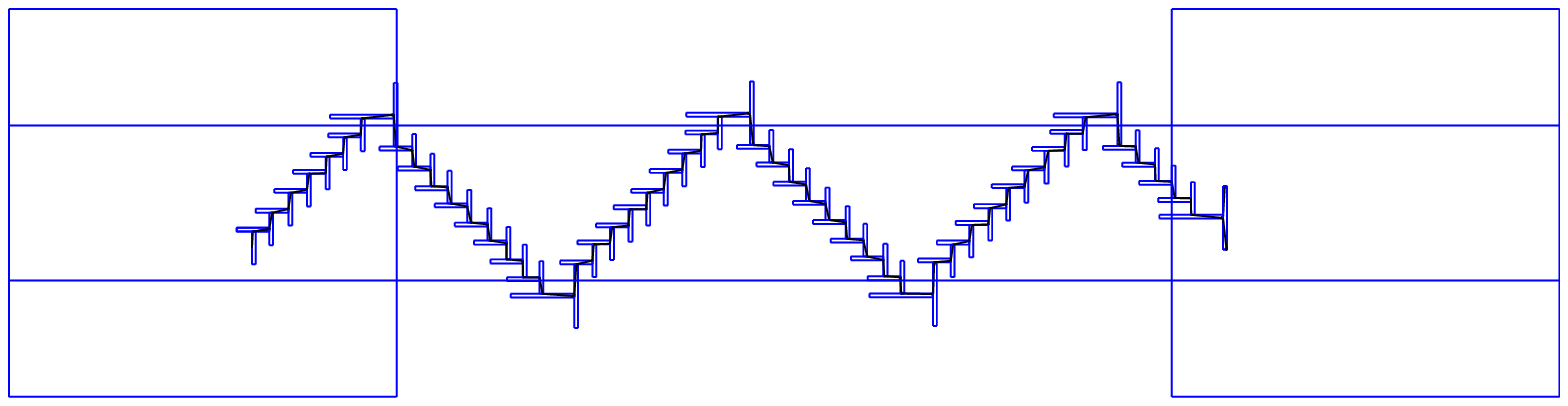}
 \caption{Producing new seeds in $R_{1,0}(x)$}
 \end{figure}

Next, we iterate the above procedure many times in both directions (to the right and to above). See Figure 9 for intuition. For any $n\in\NN$, we can construct a route from this iteration in order to get some $y,z\in R_{n,n}$ and $t,u<\infty$ through the route, such that the seed $(x\times s)_r$ is joined to the seeds $(y\times t)_r$ and
$(z\times u)_r$ within $\lceil a, b+n\MM h(1+\ii)\rfloor$.

 \begin{figure}[H]\label{picc:9}
 \center
 \includegraphics[width=13.0true cm, height=14true cm]{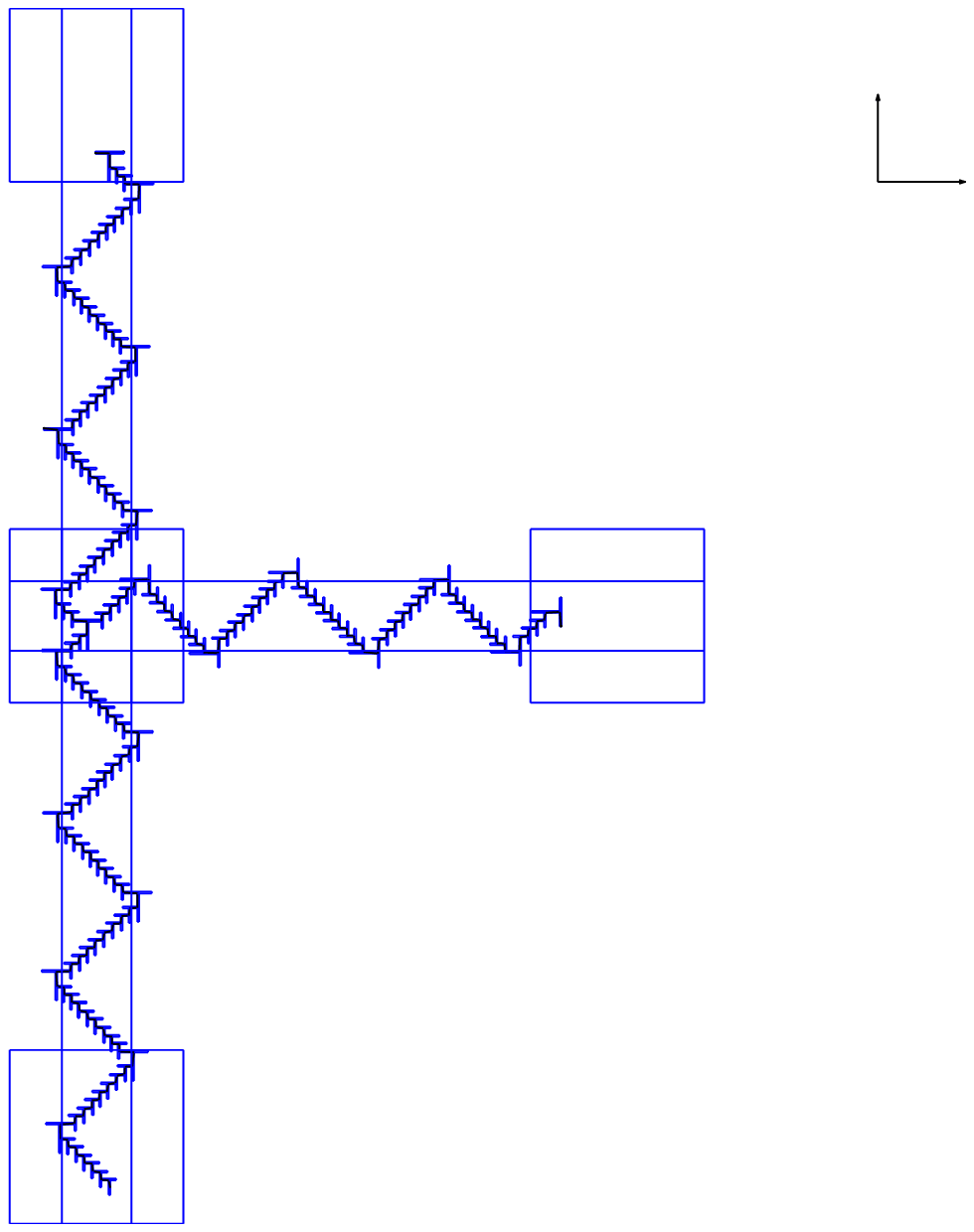}
 \caption{All S-boxes and L-boxes are disjoint}
 \end{figure}

For any valid sample (that is, a route can be successfully found), we can let the route be unique in some manner. For example, if both the seed in $R_{i-1,j}(x)$ and the seed in $R_{i,j-1}(x)$ can generate new seeds in $R_{i,j}(x)$ in finite time, then we choose the route from $R_{i-1,j}(x)$ to $R_{i,j}(x)$. That is, we put priority to the \lq left neighbor\rq. See Figure 10 for intuition. From this, we can get that there exist $y,z\in
R_{n,n}$, such that the seed $(x\times s)_r$ is joined to two seeds
$(y\times t_1^{(n)})_r$ and $(z\times t_2^{(n)})_r$ within $\lceil
a, b+n\MM h(1+\ii)\rfloor$. Furthermore, $t_1^{(n)}+t_2^{(n)}<\infty$ with large probability (depending on $n$). Denote
$$F_1(s, x, n,1+\ii):=t_1^{(n)}~~\text{and~~}
F_2(s, x, n,1+\ii):=t_2^{(n)},$$ where $1+\ii$ indicates that
the orientation of infection is northeast.

 \begin{figure}[H]\label{picc:10}
 \center
 \includegraphics[width=10.0true cm, height=10true cm]{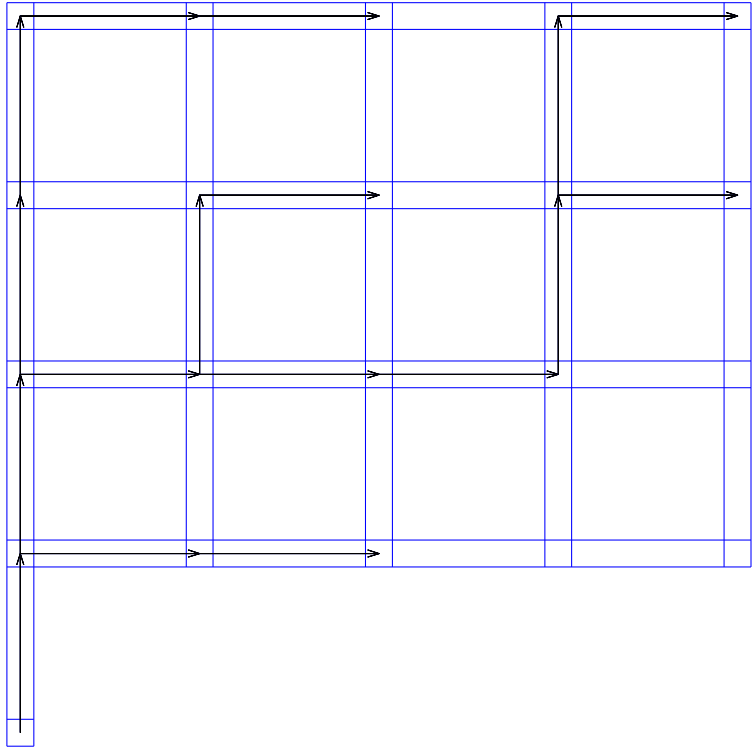}
 \caption{Dynamic renormalization ($n=4$)}
 \end{figure}

Similarly, we can define $F_1(s, x, n,o)$ and $F_2(s, x,
n,o)$ for other orientations $o\in\{1-\ii,-1+\ii,-1-\ii\}$. If $F_1(s, x, n,o)+F_2(s, x, n,o)<\infty$,
then there exist $x_1,x_2\in \lceil a, b\rfloor+n\MM ho$, such
that the seed $(x\times s)_r$ is joined to two seeds $(x_1\times
F_1(s, x, n,o))_r$ and $(x_2\times F_2(s, x, n,o))_r$, and $x, x_1,x_2$ are arranged {\it
clockwise}.\\

Having made the above preparations, we can now state the main
proposition in this section.
\begin{prop}\label{p:4.1}Suppose that $\P(\xi^0\text{ survives})>0$. Let
$x=x(\varepsilon)\in \HH$ with $\Im(x)>10h$, and let $(x\times 0)_r$ be a
horizontal seed. Then there exists  $\overline{W}>0$ which depends
only on $\varepsilon$ and $\lambda$, such that
$$\lim_{\varepsilon\rightarrow 0+}\liminf_{n\rightarrow\infty}\P\left( \frac{7\overline{W}}{6}n<F_1(0, x,
n,1+\ii)<\frac{11\overline{W}}{6}n\right)=1$$ and
$$\lim_{\varepsilon\rightarrow 0+}\liminf_{n\rightarrow\infty}\P\left( \frac{7\overline{W}}{6}n<F_2(0, x,
n,1+\ii)<\frac{11\overline{W}}{6}n\right)=1,$$
where $F_1(0,x,n,1+\ii)$ and $F_2(0,x,n,1+\ii)$ are the time points that generate the two seeds in $R_{n,n}(x)$ from the original seed $(x\times0)_r$, respectively, as defined above.
\end{prop}

The proof of Proposition \ref {p:4.1} is quite similar to the proof of Proposition 4.1 in Chen and Yao \cite{Chen-Yao2009}. So we omit the formal proof here. Readers can refer to Appendix 2 in Chen and Yao \cite{Chen-Yao2009} for details. We only
state the idea here. We have got a route by which a seed in
$R_{m,n}(x)$ is joined to other seeds in $R_{m+1,n}(x)$ and $R_{m,n+1}(x)$
with large probability. As a result, we use the \lq dynamic renormalization\rq~
method and consider each $R_{m,n}(x)$ as one site. Declare $R_{0,0}(x)$
open if $x\in R_{0,0}(x)$ and $(x\times 0)_r$ is a seed. For
$m+n\geq1$, declare $R_{m,n}(x)$ open if and only if one of the following holds.
\newline (i) $R_{m-1,n}(x)$ is open and the seed in $R_{m-1,n}(x)$ is
joined to two seeds in $R_{m,n}(x)$.
\newline(ii) $R_{m-1,n}(x)$ is closed,  $R_{m,n-1}(x)$ is  open, and the seed in $R_{m,n-1}(x)$ is joined to two seeds in $R_{m,n}(x)$.
\newline Refer to Figure 10 for intuition. The process
$(R_{m,n}(x))_{m\in\ZZ,~n\in\ZZ^+}$ is thus an oriented site percolation. Refer to
Durrett \cite{Durrett1984} and Grimmett \cite{Grimmett1999} for more detailed introductions. We can then find a unique open path from $R_{0,0}(x)$ to
$R_{n,n}(x)$ with large probability. Furthermore, we can find the
unique route constructed by S-boxes and L-boxes, within which the
seed in $R_{0,0}(x)$ is joined to another two seeds in $R_{n,n}(x)$. This
implies that $F_1(s, x, n,1+\ii)$ is the sum of the times spent in each
box. And $F_2(s, x, n,1+\ii)$ also. Figure 9 indicates that all S-boxes and L-boxes are disjoint.
So the times spent in each box are independent under certain
conditions (this has been used several times in Section 3; for details, readers can refer to Lemmas 3.5 and 3.6 of \cite{Chen-Yao2009}). Through rigorous calculation, we
get that the total number of S-boxes on the route is between
$2nj_{lower}$ and $2nj_{upper}$. Then, by the law of large numbers,
with large probability, the time spent in these S-boxes is between
$\frac{7}{6}\overline{S}n$ and $\frac{11}{6}\overline{S}n$. We can
deduce that with large probability, the time spent in these L-boxes is between
$\frac{7}{6}\overline{L}n$ and $\frac{11}{6}\overline{L}n$, too.
Hence with large probability, the total time $F_1(s, x, n,1+\ii)$ is between
$\frac{7}{6}\overline{W}n$ and $\frac{11}{6}\overline{W}n$. And $F_2(s, x, n,1+\ii)$ also. Here
$j_{lower}$ and $j_{upper}$ are two constants which satisfy $1\le
j_{upper}/j_{lower}<\frac{11}{6}$,
  and $\overline{S}$, $\overline{L}$ and $\overline{W}$
 depend only on $\lambda$ and $\varepsilon$.

\section{The complete convergence theorem}
Having established the dynamic renormalization construction, we are now in a position to prove the complete convergence theorem, Theorem \ref{t:1.1}. By Theorem 1.12 of Liggett \cite{Liggett1999}, to prove Theorem \ref{t:1.1} it suffices to prove that there exists $\Omega_0\subseteq\Omega_1$ with
$\P^{\mu}(\Omega_0)=1$, such that, for all $\omega\in\Omega_0$, the next two assertions hold.

(a)~$\P_{\lambda}\left(x\in\limsup\limits_{t\rightarrow\infty}\xi_t^A(\lambda)\right)=\P_{\lambda}(\xi^A(\lambda)\text{ survives})$ for all $x\in \HH$ and $A\subset\HH$.

(b)~$\lim\limits_{l\rightarrow\infty}\liminf\limits_{t\rightarrow\infty}\P_{\lambda}(\xi_t^{B_x(l)}(\lambda)\cap
B_x(l)\neq\emptyset)=1$ for all $x\in\HH$.\\

We will prove (a) and (b) rigorously in Sections 5.1 and 5.2, respectively. The intuitive idea is as follows. We iterate the construction posed in Proposition \ref{p:4.1} four times to get that, with large probability, a seed in $\lceil a,b\rfloor\times0$ is joined to another seed in $\lceil e,f\rfloor\times[3\overline{W}n,\infty)$. See Figure 11 for intuition. From this, we get (a). Extra tricks are needed to check (b). We will prove that, for each $n$, with large probability, a seed in $\lceil\widetilde{a}, \widetilde{b}\rfloor\times[0,\widetilde{W}]$ is joined  to another seed in $\lceil\widetilde{e},\widetilde{f}\rfloor\times[(n-1)\widetilde{W},(n+1)\widetilde{W}]$.
Together with the fact  that every remote site  cannot be infected
in a short time, we get (b).

\begin{figure}[H]\label{picc:11}
\center
\includegraphics[width=10.0true cm, height=10true cm]{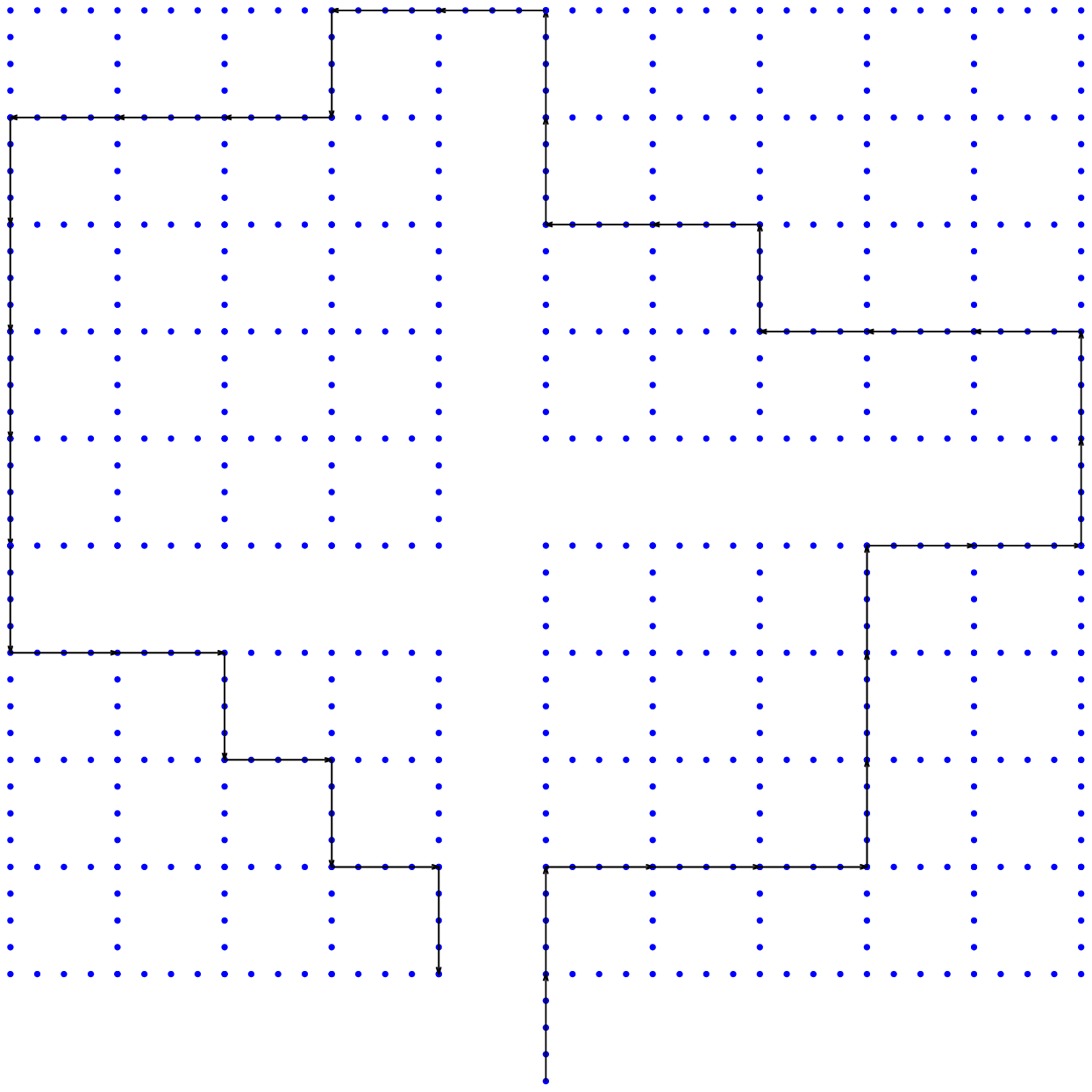}
\caption{Description of (a) ($m=5$)}
\end{figure}

\subsection{Proof of (a)}
Without loss of generality, we suppose that $\P(\xi^A\text{ survives})>0$, since otherwise both sides in (a) are equal to $0$ and (a) holds trivially. We first prove the case when $A$ is a nonempty \emph{finite} subset
of $\HH$. Let $x_0$ be any element of $A$, and let $\sigma_0=0$. Hence
$x_0$ is infected at time $\sigma_0$ for the process $\xi^A$.
Then define $\delta_k$, $\tau_k$, $Y_k$, $\sigma_{k+1}$, and $x_{k+1}$
inductively for $k\ge 0$ as follows. (See Figure 12 for intuition.) Let
$$\delta_k:=\sup\{t>\sigma_k:~x_k\times\sigma_k\text{ is joined within }\ang{x_k-r-1,
x_k+r+1+2000h\ii}$$
$$\text{to }\lceil x_k-r-1,
x_k+r+1+2000h\ii\rfloor\times t\}$$
be the death time for the contact process starting with single infection $x_k$ at time $\sigma_k$ and evolving within $\ang{x_k-r-1,
x_k+r+1+2000h\ii}$. Then $\delta_k<\infty$ almost surely on $\{\sigma_k<\infty\}$. Let
$$\tau_k:=\min\{t>\sigma_k:~x_k\times \sigma_k\text{ is joined
within }\ang{x_k-r-1, x_k+r+1+2000h\ii}$$
$$\text{to }z\times t\text{ for all }z\in \lceil x_k-r+2000h\ii,x_k+r+2000h\ii\rfloor\}-\sigma_k$$
be the waiting time until the first seed  on the top appears. Let
$$Y_k:=\sup\left\{\Im(x):~x\in\bigcup\limits_{t\le \delta_k}  \xi^A_{t}\right\}.$$ Then
$Y_k<\infty$ almost surely on $\{\sigma_k<\infty\}$. Furthermore, let
 $$\sigma_{k+1}:=\inf\{t>\delta_k:~\exists x\in \xi^A_t,\text{ s.t. }\Im(x)=Y_k+2\},$$ and let $x_{k+1}$ be the corresponding infected site. Note that, for any $k=0,1,2,\cdots$, if $\tau_k<\infty$, then $\sigma_k+\tau_k<\delta_k$.
 \begin{figure}[H]\label{picc:15}
 \center
 \includegraphics[width=15.0true cm, height=10true cm]{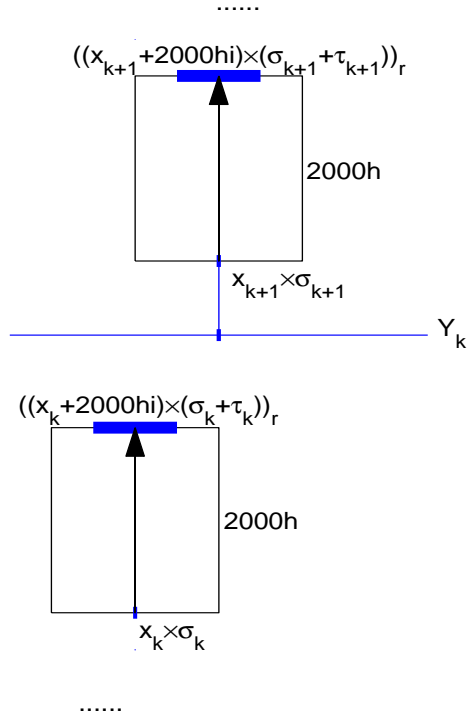}
 \caption{Inductive definitions}
 \end{figure}

Define $K:=\min\{k:~\tau_k<\infty\}$ and denote $p:=\P(\tau_0<\infty)>0$. For $t>0$, we use $\AA_t$ to denote the $\sigma$--fields generated by the graphical representation for the contact process until time $t$. Therefore, by translation invariance and the fact that $\sigma_k$ is a stopping time for all $k\in\NN$, we get that, if $\sigma_k<\infty$ for all $k$ and $\sigma_k\uparrow\infty$, then
$$\P(K<\infty~|~\AA_{\sigma_k})\geq\P(\tau_k<\infty~|~\AA_{\sigma_k})=\P(\tau_0<\infty)=p>0$$
for any $k=0,1,2,\cdots$. That is,
\begin{equation}\label{e:6.0}
\P(\P(K<\infty~|~\AA_{\sigma_k})\geq p>0\text{ for all }k=0,1,2,\cdots~|~\sigma_k<\infty\text{ for all }k,~\sigma_k\uparrow\infty)=1.
\end{equation}
Furthermore, using the martingale convergence theorem, we can get that
$$\P(K<\infty~|~\AA_{\sigma_k})\rightarrow\textbf{1}_{\{K<\infty\}}\text{ a.s. on }\{\sigma_k<\infty\text{ for all }k,~\sigma_k\uparrow\infty\}$$
as $k$ tends to infinity. So, by (\ref{e:6.0}), we get
\begin{equation}\label{e:6.01}
\P(K<\infty~|~\sigma_k<\infty\text{ for all }k,~\sigma_k\uparrow\infty)=1.
\end{equation}
Also, note that
\begin{equation}\label{e:6.02}
\P(\sigma_k<\infty\text{ for all }k,~\sigma_k\uparrow\infty~|~\xi^A\text{ survives})=1.
\end{equation}
By (\ref{e:6.01}) and (\ref{e:6.02}), together with our assumption that $\P(\xi^A\text{ survives})>0$, we get
\begin{equation}\label{e:7.1}
\P(K<\infty~|~\xi^A\text{ survives})=1. \end{equation}

If $K<\infty$, then let $y_1:=x_K+2000h\ii$, and let $t_1=\sigma_k+\tau_k$. Therefore, $(y_1\times t_1)_r$ is a horizontal seed. Let
$$\zeta=F_1(F_1(F_1(F_1(t_1,y_1,m,1+\ii),y_2,m,-1+\ii),y_3,m,-1-\ii),y_4,m-1,1-\ii),$$
and let $(\vartheta\times \zeta)_r$ be the corresponding seed if
$\zeta<\infty$. Here, $F_1$ is defined as in Section 4, and $y_2$, $y_3$, $y_4$ are the centers of corresponding seeds in each step. Therefore, $$\vartheta\in R_{-1,1}(y_1)\subset
B_{y_1}(2\MM h).$$ See Figure 11 for intuition. Note that $\zeta$ is the sum of the times spent in each of the four orientations as shown in Figure 11. These times are independent under certain conditions (this has been used several times in Section 3; for details readers can refer to Lemmas 3.5 and 3.6 of \cite{Chen-Yao2009}). Together with Proposition \ref{p:4.1}, we get
\begin{align*}
\lim_{\varepsilon\rightarrow
0+}\liminf_{m\rightarrow\infty}\P(3\overline{W}m\le
\zeta<\infty~|~K<\infty)=1,
\end{align*}
which implies that
$$
\lim_{\varepsilon\rightarrow
0+}\liminf_{m\rightarrow\infty}\P(\exists t\ge
3\overline{W}m ,~\text{s.t. }\xi^A_t\cap B_{y_1}(2\MM h)\neq
\emptyset~|~K<\infty)=1.
$$
By the dominated convergence theorem, we have
$$
\lim_{\varepsilon\rightarrow
0+}\P\left(\left.\limsup_{t\rightarrow\infty} \xi^A_t\cap
B_{y_1}(2\MM h)\neq \emptyset~\right|~K<\infty\right)=1.
$$
Furthermore,
\begin{equation}\label{e:7.2}
\lim_{\varepsilon\rightarrow
0+}\P\left(\left.\limsup_{t\rightarrow\infty} \xi^A_t\neq
\emptyset~\right|~K<\infty\right)=1.
\end{equation}
By (\ref{e:7.1}) and (\ref{e:7.2}), we have
$$
\P\left(\limsup_{t\rightarrow\infty} \xi^A_t\neq
\emptyset~|~\xi^A \text{~survives~}\right)=1.
$$
Turning to the quenched law, there exists $\Omega_A\subseteq\Omega_1$ with
$\P^{\mu}(\Omega_A)=1$, such that, for all $\omega\in \Omega_A$,
\begin{equation}\label{e:7.3}
\P_\lambda\left(\xi^A(\lambda) \text{ survives},~\limsup\limits_{t\rightarrow\infty}\xi_t^A(\lambda)=\emptyset\right)=0.
\end{equation}
That is, $\xi^A(\lambda)$ survives strongly if it survives. See
page 42 of Liggett \cite{Liggett1999} for the definition of \lq strong
survival\rq.

Fix $\omega\in\Omega_A$. For any $y,z\in\HH$, we have
$$
\P_\lambda(z\in \xi_1^y(\lambda))>0.
$$
We can construct an appropriate sequence of stopping times and use the strong Markov property under the quenched law to get
$$z\in\limsup\limits_{t\rightarrow\infty}\xi_t^A(\lambda)\text{ a.s. on }\{y\in\limsup\limits_{t\rightarrow\infty}\xi_t^A(\lambda)\}.$$
That is,
$$\P_\lambda\left(z\not\in\limsup\limits_{t\rightarrow\infty}\xi_t^A(\lambda),~y\in\limsup\limits_{t\rightarrow\infty}\xi_t^A(\lambda)\right)=0$$
for any $y,z\in\HH$. Since
$$\left\{\limsup\limits_{t\rightarrow\infty}\xi_t^A(\lambda)\neq\emptyset\right\}=\bigcup\limits_{y\in\HH}\left\{y\in\limsup\limits_{t\rightarrow\infty}\xi_t^A(\lambda)\right\},$$ we have
\begin{equation}\label{e:7.4}
\P_\lambda\left(z\not\in\limsup\limits_{t\rightarrow\infty}\xi_t^A(\lambda),~\limsup\limits_{t\rightarrow\infty}\xi_t^A(\lambda)\neq\emptyset\right)=0.
\end{equation}
From (\ref{e:7.3}) and (\ref{e:7.4}), together with the fact that $\left\{z\in\limsup\limits_{t\rightarrow\infty}\xi_t^A(\lambda)\right\}\subseteq\left\{\xi^A(\lambda)\text{ survives}\right\}$, we can deduce that, for any finite subset
  $A\subseteq\HH$, $\omega\in\Omega_A$, and
 $z\in\HH$,
$$
\P_\lambda\left(z\in \limsup_{t\rightarrow\infty}\xi_t^A(\lambda)\right)=\P\left(\xi^A(\lambda) \text{~survives}\right).
$$

Then, let $$\Omega_0':=\bigcap\limits_{A\subset
\HH,~|A|<\infty}\Omega_A.$$ Then $\P^{\mu}(\Omega_0')=1$. Moreover,
(a) holds for all $\omega\in \Omega_0'$, $x\in\HH$, and $A\subset
\HH$ with $|A|<\infty$.\\

Next, we consider the case when $|A|=\infty$. We can get that, for any $n>0$, there exists $m_n$
such that $\P(\xi^{B}\text{ survives})>1-4^{-n}$ for any $B\subset \HH$ with $|B|\ge m_n$,
for a reason similar to the proof of Lemma \ref{l:3.1}. This implies that
$$
\P^{\mu}(\{\omega\in\Omega_1:~\P_\lambda(\xi^B(\lambda)\text{ survives})\ge
1-2^{-n}\})\ge 1- 2^{-n}.
$$
Let $$\Xi_n':=\{\omega\in\Omega_1:~\P_\lambda(\xi^B(\lambda)\text{ survives})\ge
1-2^{-n}\}$$ for any $n\in\NN$. Then $\Xi_n'$ decreases as $n$ increases. Set
$$\Omega_0'':=\Omega_0'\cap\left(\bigcap\limits_n\Xi_n'\right).$$ Then $\P^{\mu}(\Omega_0'')=1$. If
$\omega\in\Omega_0''$, $x\in\HH$, $A\subset\HH$, and $|A|=\infty$, then
let $(A_n)$ be an increasing sequence of finite sets which satisfy
$\lim\limits_{n\rightarrow\infty} A_n=A$ and $|A_n|>m_n$ for all $n$. Then, for any $x\in\HH$, we have
\begin{align*}
\P_\lambda\left(x\in
\limsup_{t\rightarrow\infty}\xi^A_t(\lambda)\right)&\ge
\lim_{n\rightarrow\infty} \P_\lambda\left(x\in
\limsup_{t\rightarrow\infty}\xi^{A_n}_t(\lambda)\right)\\&=\lim_{n\rightarrow\infty}
\P_\lambda\left( \xi^{A_n}_t(\lambda) \text {~survives}\right)\\&\ge
\lim_{n\rightarrow\infty} (1-2^{-n})=1.
\end{align*}
But $\xi^A(\lambda)$ survives with $\P_{\lambda}$-probability $1$. As a result,
$$\P_{\lambda}(x\in\limsup\limits_{t\rightarrow\infty}\xi^{A}_t(\lambda))=\P_\lambda(\xi^{A}_t(\lambda)\text{ survives})=1.$$
Furthermore, (a) holds for all $\omega\in\Omega_0''$, $x\in\HH$, and $A\subset\HH$.

\subsection{Proof of (b)}
We begin with the seed $(x\times
s)_r$. For convenience, for any $n\in\NN$, we use the following algorithm to generate a new seed from $(x\times s)_r$ and record the time used. Recall that, in the algorithm, $F_1$ and $F_2$ are as defined in Section 4.

\noindent\textbf{Algorithm}\\
0) Set $t=s$ and $y=x$.\\
1) Set $s'=s-100\overline{W}n[s/100\overline{W}n], ~  v=8\cdot
\textbf{1}_{\{s'\le 37\overline{W}n\}}$ and $u=9-v$.\\
 One can check that
   $$s'+(6u+10v)\cdot
   \left[\frac{7}{6}\overline{W}n,\frac{11}{6}\overline{W}n\right]\subseteq[100\overline{W}n,200\overline{W}n).$$
Operate 2)$\sim$7) $u$ times\\
2) $t=F_2(t,n,1+\ii)$;\\
3) $t=F_1(t,n,1-\ii)$;\\
4) $t=F_1(t,n,1+\ii)$;\\
5) $t=F_1(t,n,-1+\ii)$;\\
6) $t=F_2(t,n-1,-1-\ii)$;\\
7) $t=F_2(t,n+1, -1+\ii)$;\\
Operate 8)$\sim $17) $v$ times\\
8) $t=F_2(t,n,1+\ii)$;\\
9) $t=F_1(t,n,1-\ii)$;\\
10) $t=F_2(t,n,1+\ii)$;\\
11) $t=F_1(t,n,1-\ii)$;\\
12) $t=F_1(t,n,1+\ii)$;\\
13) $t=F_1(t,n, -1+\ii)$;\\
14) $t=F_2(t,n-1,-1-\ii)$;\\
15) $t=F_1(t,n,-1+\ii)$;\\
16) $t=F_2(t,n, -1-\ii)$;\\
17) $t=F_2(t,n+1,-1+\ii)$;\\
18) Return $t$.\\

If the output value $t<\infty$, then the corresponding site belongs to
$R_{18(n+1),0}(x)$. Moreover, by Proposition \ref{p:4.1}, we know that  $t\in[100\overline{W}n,200\overline{W}n)$ with large probability if $s\in
[0,100\overline{W}n)$. Denote
$$G(s,x,n,\ii):=t.$$
See Figure 13 for intuition. Similarly, we denote
$G(s,x,n,1)$ the corresponding site that belongs to
$R_{0,18(n+1)}(x)$ generated in the same way (but in a different direction).\\

\begin{figure}[H]\label{picc:12}
\center
\includegraphics[width=15.0true cm, height=10true cm]{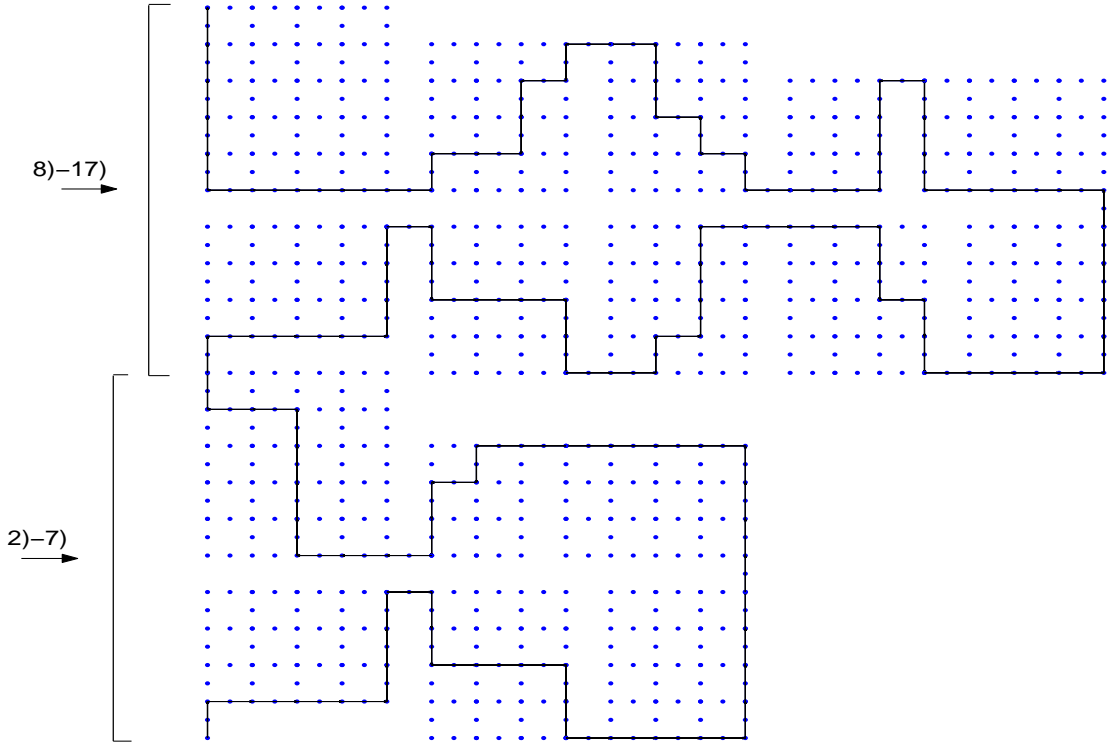}
\caption{Description of $G(s,x,n,\ii)$}
\end{figure}

Next, we iterate the above procedure many times in both directions (to the right and above).
See Figure 10 for intuition. For any $m\in\NN$, we can construct a route through this iteration in order to get a new seed in $R_{18(n+1)m,18(n+1)m}(x)$.
The procedure is similar to the argument before Proposition \ref{p:4.1}, and we can use a similar way (prior to the \lq left neighbor\rq) to make the route unique. We denote the time by $L(s,x,n,m,1+\ii)$, which is finite with large probability (depending on $n$ and $m$). Here, $1+\ii$ indicates that the orientation of infection is northeast.
%
%
%
%

 Similarly, we can define
 $L(s,x,n,m,o)$ for other orientations $o\in\{1-\ii,-1+\ii,-1-\ii\}$. We then have the following proposition, which is parallel to Proposition \ref{p:4.1}, but it is more accurate.
\begin{prop}\label{p:5.1}Suppose that $\P(\xi^0\text{ survives})>0$. Let $x=x(\varepsilon)\in \HH$ with
$\Im(x)>10h$, and let $(x\times 0)_r$ be a horizontal seed. Then
$$\lim_{\varepsilon\rightarrow 0+}\liminf_{n\rightarrow\infty}\liminf_{m\rightarrow\infty}
\P(200\overline{W}nm<L(0,x,n,m,1+\ii)<200\overline{W}n(m+1))=1.$$
\end{prop}
\pf By Proposition \ref{p:4.1} and the FKG inequality, we have that with large probability
$$G(s,x,n,\ii)\in[100k\overline{W}n,100(k+1)\overline{W}n)~~\text{and}~~G(s,x,n,\ii)\in[100k\overline{W}n,100(k+1)\overline{W}n)$$
if $s\in [100(k-1)\overline{W}n,100k\overline{W}n)$. Similar to the idea of Proposition \ref{p:4.1}, this situation corresponds to a 1-dependent site percolation. Using the result of 1-dependent site percolation (see \cite{Durrett1984}), we get the conclusion.\qed\\

Next, we prove (b). Without loss of generality, we suppose that $\Im(x)\ge 10h$. Suppose that $(x\times0)_{r}$ is a horizontal seed. Let
$$\mu:=L(L(L(L(0,x,n,m,1+\ii),x_1,n,m,-1+\ii),x_2,n,m,-1-\ii),x_3,n,m-1,1-\ii),$$
and let $(\nu\times\mu)_r$ be the corresponding seed if $\mu<\infty$. Then $\nu\in B_x(40n\MM h)$. Here $x_1$, $x_2$, and $x_3$ are the centers of the corresponding seeds in each step. Note that $\mu$ is the sum of the times spent in each of the four orientations. These times are independent under certain conditions (this has been used several times in Section 3; for details readers can refer to Lemmas 3.5 and 3.6 of \cite{Chen-Yao2009}). Together with Proposition \ref{p:5.1}, we get
$$
\lim_{\varepsilon\rightarrow
0+}\liminf_{n\rightarrow\infty}\liminf_{m\rightarrow\infty}\P(800\overline{W}nm-200\overline{W}n\le\mu
\le800\overline{W}nm+600\overline{W}n)=1.
$$
That is,
$$
\lim_{\varepsilon\rightarrow
0+}\liminf_{n\rightarrow\infty}\liminf_{m\rightarrow\infty}\P
(\exists t\in [800\overline{W}n(m-1), 800\overline{W}n(m+1)],~\text{s.t. }
\xi_t^{\lceil x-r,x+r\rfloor}\cap B_x(40n\MM h)\neq
\emptyset)=1.
$$
We can deduce that, for any $\delta>0$, there exist $n_0\in\NN$ and
$m_0\ge 2$, such that
$$
\P(\exists t\in [800\overline{W}n_0(m_0-1),
800\overline{W}n_0(m_0+1)],\text{ s.t. }\xi_t^{\lceil x-r,x+r\rfloor}\cap
B_x(40n_0\MM h)\neq \emptyset)>1-\delta^2.
$$
Turning to the quenched law, denote
$$\Omega_\delta^{(1)}:=\{\omega\in\Omega_1:~\P_\lambda(\exists t\in [800\overline{W}n_0(m_0-1),
800\overline{W}n_0(m_0+1)],\text{ s.t. }\xi_t^{\lceil x-r,x+r\rfloor}\cap
B_x(40n_0\MM h))>1-\delta\}.$$ Then $\P^\mu(\Omega_\delta^{(1)})\geq1-\delta$.

On the other hand, consider the Richardson's process $(\zeta_t)$ on $\HH$ by suppressing all recoveries from $(\xi_t)$, we have
$$\lim\limits_{l\rightarrow\infty}\P(\inf\{t>0:~\zeta_t^A\cap B_x(40n_0\MM h)\neq\emptyset\}\leq800\overline{W}n_0(m_0+1)+1\text{ for some finite subset }A\subseteq\HH\setminus B_x(l))=0.$$
So for the above $\delta>0$, there exists $l_\delta>40n_0\MM h$, such that
$$\P(\text{for some finite subset }A\subseteq\HH\setminus B_x(l_\delta),~\text{there exists }t\in(0,800\overline{W}n_0(m_0+1)+1]$$
$$\text{s.t. }\zeta_t^A\cap B_x(40n_0\MM h)\neq\emptyset)<\delta^2.$$
Turning to the quenched law, denote
$$\Omega_\delta^{(2)}:=\{\omega\in\Omega_1:~\P_\lambda(\text{for some finite subset }A\subseteq\HH\setminus B_x(l_\delta),~\text{there exists }t\in(0,800\overline{W}n_0(m_0+1)+1]$$
$$\text{s.t. }\zeta_t^A\cap B_x(40n_0\MM h)\neq\emptyset)<\delta\}.$$
Then $\P^\mu(\Omega_\delta^{(2)})\geq1-\delta$. So $\P^\mu(\Omega_\delta^{(1)}\cap\Omega_\delta^{(2)})\geq1-2\delta$.

Next, fix $\omega\in\Omega_\delta^{(1)}\cap\Omega_\delta^{(2)}$. For any $s\geq1$, set
$$
\tau_s:=\inf\{u\ge s-1:~\xi_u^{\lceil x-r,x+r\rfloor}\cap
B_x(l_\delta)=\emptyset\}.
$$
Then $\tau_s$ is a stopping time. Using the strong Markov property under the quenched law, together with the facts that $\xi_t^A\subseteq\zeta_t^A$ for any $t$ and $\zeta_t^A$ increases as $t$ increases, we can get that, for any finite subset $A\subseteq\HH\setminus B_x(l_\delta)$,
\begin{align*}
&\P_{\lambda}(\exists 0<t\leq800\overline{W}n_0(m_0+1),\text{ s.t. }\xi_{t+s}^{\lceil
x-r,x+r\rfloor}(\lambda)\cap B_x(40n_0\MM h)\neq
\emptyset~|~\xi_{\tau_s}^{\lceil x-r,x+r\rfloor}(\lambda)=A )\\
\le&\P_\lambda(\exists 0<t\leq800\overline{W}n_0(m_0+1)+1,\text{ s.t. }\zeta_{t}^{A}\cap
B_x(40n_0\MM h)\neq \emptyset)\le\delta.
\end{align*}
Then we use the strong Markov property under the quenched law again to get
\begin{align*} &\P_{\lambda}(\exists
0<t\leq800\overline{W}n_0(m_0+1),\text{ s.t. }\xi_{t+s}^{\lceil
x-r,x+r\rfloor}(\lambda)\cap B_x(40n_0\MM h)\neq
\emptyset,~\xi_u^{\lceil x-r,x+r\rfloor}\cap B_0(l_\delta)= \emptyset\\
&~~~~~~~~~~~~~~~~~~~~~~~~~~~~~~~~~~~~~~~~~~~~~~~~~~~~~~~~~~~~~~~~~~~~~~~~~~~~~~~~~~~~~~~~~~\text{ for some }u\in [s-1,
s])\\
=&\P_\lambda(\exists 0<t\leq800\overline{W}n_0(m_0+1),\text{ s.t. }\xi_{t+s}^{\lceil
x-r,x+r\rfloor}(\lambda)\cap
B_x(40n_0\MM h)\neq \emptyset,~\tau_s\le s)\\
=&\P_\lambda(\P_\lambda(\exists
0<t\leq800\overline{W}n_0(m_0+1),\text{ s.t. }\xi_t^{\lceil x-r,x+r\rfloor}(\lambda)\cap
B_x(40n_0\MM h)\neq \emptyset~|~\FF_{\tau_s} );~\tau_s\le s)\\
=&\P_\lambda(\P_\lambda(\exists
0<t\leq800\overline{W}n_0(m_0+1),\text{ s.t. }\xi_t^{\lceil x-r,x+r\rfloor}(\lambda)\cap
B_x(40n_0\MM h)\neq \emptyset~|~\xi_{\tau_s}^{\lceil x-r,x+r\rfloor}(\lambda) );~\tau_s\le s)\\
\le &\delta\cdot\P_\lambda(\tau_s\le s)\leq\delta
\end{align*} for any
$s\ge 1$. Therefore, for any $s\geq 1$, we have
\begin{align*}
&\P_\lambda(\xi_u^{\lceil x-r,x+r\rfloor}\cap B_x(l_\delta)\neq
\emptyset \text{ for all }u\in [s-1, s])\\
\geq&\P_\lambda(\exists
0<t\leq800\overline{W}n_0(m_0+1),\text{ s.t. }\xi_{t+s}^{\lceil x-r,x+r\rfloor}(\lambda)\cap
B_x(40n_0\MM h)\neq \emptyset)\\
&-\P_\lambda(\exists
0<t\leq800\overline{W}n_0(m_0+1),\text{ s.t. }\xi_{t+s}^{\lceil x-r,x+r\rfloor}(\lambda)\cap
B_x(40n_0\MM h)\neq \emptyset,~\xi_u^{\lceil x-r,x+r\rfloor}\cap
B_x(l_\delta)=
\emptyset\\ &~~~~~~~~~~~~~~~~~~~~~~~~~~~~~~~~~~~~~~~~~~~~~~~~~~~~~~~~~~~~~~~~~~~~~~~~~~~~~~~~~~~~~~\text{ for some }u\in [s-1, s] )\\
\geq&1-2\delta.
\end{align*}
Since $\xi_t^{A_1}\subseteq \xi_t^{A_2}$ for all $t\geq0$ if $A_1\subseteq A_2$, we
have, for any $\omega\in\Omega_\delta^{(1)}\cap\Omega_\delta^{(2)}$,
$$
\P_\lambda(\xi_s^{B_x(l_\delta)}\cap B_x(l_\delta)\neq \emptyset\text{ for all }s\in[t,t+1])\ge 1-2\delta
$$
for any $t\geq0$. So, if we denote
$$\Omega_\delta:=\{\omega\in\Omega_1:~\P_\lambda(\xi_s^{B_x(l_\delta)}\cap B_x(l_\delta)\neq \emptyset\text{ for all }s\in[t,t+1]\geq1-2\delta)\},$$
then $$\P^\mu(\Omega_\delta)\geq\P^\mu(\Omega_\delta^{(1)}\cap\Omega_\delta^{(2)})\geq1-2\delta.$$
And furthermore, there exists $l_n\uparrow\infty$ such that, for any $n\in\NN$ and $t\geq0$,
$$\P^\mu(\Omega_{n,t})\geq1-2^{-n-t-1},$$
where we set
$$\Omega_{n,t}:=\{\omega\in\Omega_1:~\P_\lambda(\xi_s^{B_x(l_n)}\cap B_x(l_n)\neq \emptyset\text{ for all }s\in[t,t+1]\geq1-2^{-n-t-1})\}$$
for any $n\in\NN$ and $t\geq0$. Next, set
$$\Omega_n:=\bigcap\limits_{k=0}^{\infty}\Omega_{n,k}$$
for any $n\in\NN$. Then, for any $n\in\NN$, we have $\P(\Omega_n)\geq1-2^{-n}$, and, on $\Omega_n$,
\begin{align*}
\liminf\limits_{t\rightarrow\infty}\P_\lambda(\xi_t^{B_x(l_n)}\cap B_x(l_n)\neq \emptyset)&\geq\P(\forall t\geq0,~\xi_t^{B_x(l_n)}\cap B_x(l_n)\neq \emptyset)\\
&=\P_\lambda\left(\bigcap\limits_{k=0}^{\infty}\{\xi_s^{B_x(l_n)}\cap B_x(l_n)\neq \emptyset\text{ for all }s\in[k,k+1]\}\right)\geq1-2^{-n}.
\end{align*}
Note that $\Omega_n$ increases as $n$ increases. So, if we set
$$\Omega_0''':=\bigcup\limits_{n=1}^{\infty}\Omega_n,$$
then $\P^{\mu}(\Omega_0''')=1$, and, on $\Omega_0'''$,
$$\lim\limits_{n\rightarrow\infty}\liminf\limits_{t\rightarrow\infty}\P_\lambda(\xi_t^{B_x(l_n)}\cap B_x(l_n)\neq \emptyset)=1.$$
That is, (b) holds for all $\omega\in \Omega_0'''$.\\

Finally, set $\Omega_0:=\Omega_0''\cap \Omega_0'''$. As a result,
(a) and (b) hold  for all $\omega\in \Omega_0$. So, we have proved the complete convergence theorem, Theorem \ref{t:1.1}.\\

\noindent\textbf{Acknowledgements.}\quad We are grateful to the anonymous referees for their careful reading and invaluable suggestions, especially for the suggestion that shortened the proof of Lemma \ref{l:3.1}.

The first author's research was partially supported by NSFC grants (No. 11126236 and No. 10901008) and an innovation grant from ECNU. The second author's research was partially supported by NSFC grants (No. 11001173 and No. 11171218).

\end{document}